\documentclass[11pt,reqno]{article}
\usepackage{latexsym, amsmath, amssymb, a4, epsfig}
\usepackage{graphicx}
\usepackage{color}
\usepackage[swedish,english]{babel}

\usepackage{subcaption}

\newtheorem{theorem}{Theorem}[section]
\newtheorem{cor}[theorem]{Corollary}
\newtheorem{lemma}[theorem]{Lemma}

\newcommand{\nm}{\noalign{\smallskip}}
\newcommand{\qed}{ $\Box$}
\newcommand{\ds}{\displaystyle}
\newcommand{\pf}{\noindent {\sl Proof}. \ }
\newcommand{\p}{\partial}

\newcommand{\eqnref}[1]{(\ref {#1})}

\newcommand{\Cbb}{\mathbb{C}}

\newcommand{\Rbb}{\mathbb{R}}

\newcommand{\la}{\langle}
\newcommand{\ra}{\rangle}

\newcommand{\Ecal}{\mathcal{E}}

\newcommand{\Hcal}{\mathcal{H}}

\newcommand{\Kcal}{\mathcal{K}}

\newcommand{\Scal}{\mathcal{S}}

\newcommand{\Ga}{\alpha}

\newcommand{\Gd}{\delta}
\newcommand{\Ge}{\epsilon}

\newcommand{\Gvf}{\varphi}

\newcommand{\Gl}{\lambda}

\newcommand{\Gm}{\mu}
\newcommand{\Gv}{\nu}
\newcommand{\Gp}{\pi}
\newcommand{\Gt}{\theta}

\newcommand{\Gs}{\sigma}

\newcommand{\Gj}{\tau}

\newcommand{\GD}{\Delta}

\newcommand{\GO}{\Omega}

\newcommand{\beq}{\begin{equation}}
\newcommand{\eeq}{\end{equation}}

\def\ol{\overline}

\numberwithin{equation}{section}
\numberwithin{figure}{section}

\begin{document}
\title{Classification of spectra of the Neumann--Poincar\'{e} operator on planar domains with corners by resonance\thanks{\footnotesize This work is supported by A3 Foresight Program of Korea through NRF grant NRF-2014K2A2A6000567 (to H.K and M.L), by the Korean Ministry of Science, ICT and Future Planning through NRF grant No. NRF-2013R1A1A3012931 (to M.L), and by the Swedish Research
Council under contract 621-2014-5159 (to J.H).}}

\author{Johan Helsing\thanks{Centre for Mathematical Sciences, Lund University, 221 00 Lund, Sweden (helsing@maths.lth.se)}
 \and Hyeonbae Kang\thanks{Department of Mathematics, Inha University, Incheon
22212, S. Korea (hbkang@inha.ac.kr).} \and Mikyoung Lim\thanks{Department of Mathematical Sciences,
Korea Advanced Institute of Science and Technology, Daejeon 305-701, Korea (mklim@kaist.ac.kr).}}

\date{}

\maketitle

\begin{abstract}
We study spectral properties of the Neumann--Poincar\'{e} operator on planar domains with corners with particular emphasis on existence of continuous spectrum and pure point spectrum. We show that the rate of resonance at continuous spectrum is different from that at eigenvalues, and then derive a method to distinguish continuous spectrum from eigenvalues. We perform computational experiments using the method to see whether continuous spectrum and pure point spectrum appear on domains with corners. For the computations we use a modification of the Nystr\"om method which makes it possible to construct high-order convergent
discretizations of the Neumann--Poincar\'{e} operator on domains with corners. The results of experiments show that all three possible spectra, absolutely continuous spectrum, singularly continuous spectrum, and pure point spectrum, may appear depending on domains. We also prove rigorously two properties of spectrum which are suggested by numerical experiments: symmetry of spectrum (including continuous spectrum), and existence of eigenvalues on rectangles of high aspect ratio.
\end{abstract}

\noindent{\footnotesize {\bf AMS subject classifications}. 35P05(primary),  45B05(secondary)}

\noindent{\footnotesize {\bf Key words}. Neumann--Poincar\'{e} operator, Lipschitz domain, spectrum, RCIP method, resonance}

\section{Introduction}

Let $\GO$ be a bounded simply connected domain in $\Rbb^2$ with the Lipschitz boundary. The Neumann--Poincar\'{e} (NP) operator on $\GO$ is defined by
\beq\label{NPoperator}
\Kcal_{\p\GO}^* [\Gvf] (x) := \frac{1}{2\pi} \text{p.v.} \int_{\p \GO} \frac{(x-y) \cdot n(x)}{|x-y|^2} \Gvf (y) \, d\Gs(y) \;, \quad x \in \p\GO,
\eeq
where $n(x)$ denotes the unit outward normal vector to $\p\GO$ at $x \in \p\GO$ and p.v. stands for the Cauchy principal value. Recently there is rapidly growing interest in the spectral properties of the NP operator due to its relation to plasmonics and cloaking by anomalous localized resonance: Plasmon resonance occurs at eigenvalues of the NP operator and anomalous localized resonance occurs at the accumulation point of eigenvalues, respectively (see, for example, \cite{ACKLM-ARMA-13, MFZ-PR-05, MN-PRSA-06} and references therein).

Although the NP operator $\Kcal_{\p \GO}^{*}$ is not self-adjoint with respect to the usual $L^2$-inner product unless $\GO$ is a disk or a ball \cite{Lim-IJM-01}, it can be realized as a self-adjoint operator on $H^{-1/2}(\p\GO)$ space by introducing a new inner product (see \cite{Kang-SMF-15, KPS-ARMA-07} and the next section for a brief review). Here and throughout this paper $H^s$ denotes the usual $L^2$ Sobolev space. So, the NP operator can have only three kinds of spectra: absolutely continuous spectrum, singularly continuous spectrum, and pure point spectrum (eigenvalues) \cite{RS-book-80, Yosida-book-74}.

Observe that the NP operator depends on $\GO$ in two ways: integration over $\p\GO$ and the normal vector $n(x)$. So, its spectral nature differs depending on smoothness of the domain on which it is defined. If the domain has a smooth boundary, $C^{1, \Ga}$ for some $\Ga>0$ to be precise, then the NP operator is compact on $H^{-1/2}(\p\GO)$ and its spectrum consists of eigenvalues converging to $0$. It is worth mentioning that a convergence rate of eigenvalues of the NP operator on smooth domains (among others) is obtained in a recent paper \cite{MS-arXiv}.
However, not much is known about the spectrum of the NP operator defined on the domain with corners. Bounds on the essential spectrum on curvilinear polygonal domains have been obtained in \cite{PP-JAM-14}. Recently a complete spectral resolution of the NP operator on the intersecting disk has been derived in \cite{KLY}, which in particular shows that there is only absolutely continuous spectrum, no point spectrum and no singularly continuous spectrum. It is worth mentioning that it also shows that the bound in previous mentioned paper is sharp for the intersecting disks. We mention that T. Carleson solutions of the interface problem on intersecting disks, which is closely related to the continuous spectrum (there was no notion of continuous spectrum at his time). So the paper \cite{KLY} may be regarded as a modern (and more complete) treatise of \cite{Carleman-book-16} even if the former was written without knowing existence of the latter. We also refer to a recent work \cite{pendry} where plasmon resonance on intersecting disks was studied in a numerical way.

Some natural questions arise regarding the spectrum of the NP operator on general domains with corners: Does it always have a continuous spectrum?,  no point spectrum?, and so on. The purpose of this paper is to address these questions.
On the one hand, we investigate these questions in a numerical way as the first step toward a better understanding of the spectral nature of the NP operator on domains with corners. Novelty of the computational approach of this paper may be found in two aspects. Firstly we present a way to classify spectra based on resonance. Extending the analysis of \cite{KLY} we show that the resonance at absolutely continuous spectrum is weaker than that at eigenvalues. We quantify the rate of resonance and develop a computational method to distinguish absolutely continuous spectrum from eigenvalues. Singularly continuous spectrum is more difficult to classify. However, if a strong resonance occurs at a point inside the absolutely continuous spectrum, we may infer that the point is in the singularly continuous spectrum. Another important aspect of this paper is the method of computation. For classification of spectra using resonance high precision computations are required. However, when the domain has corners, it is quite difficult to compute with high precision the NP operator. In this paper we use the Recursively Compressed Inverse Preconditioning (RCIP) method \cite{HelsOjal08}, explained in detail in \cite{Hels11, Helsing-tutorial}, which is a high precision method for solving integral equations on piecewise smooth boundaries. It is worth mentioning that this method has been successfully adapted for computation of polarizability on domains with corners \cite{HMM-NJP-11}.

Results of computational experiments of this paper reveal that absolutely continuous spectrum always appears on planar domains with corners, while pure point spectrum and singularly continuous spectrum may or may not appear depending on domains. For example, for rectangles there is a threshold $r_0$ of the aspect ratio such that if the aspect ratio is less than $r_0$ then no eigenvalue appears, and if the aspect ratio is larger than $r_0$ then eigenvalues appear. In fact, more and more eigenvalues appear as the aspect ratio increases. But no singularly continuous spectrum appears. On perturbed ellipses, singularly continuous spectrum (discrete eigenvalues embedded in absolutely continuous spectrum) appears.

On the other hand, we prove rigorously some important spectral properties of the NP operator suggested by computational experiments. We first show that the spectrum of the NP operator on planar domains is symmetric with respect to $0$. This fact is known for eigenvalues \cite{Schiffer-PJM-57}. We extend it to include the continuous spectrum by proving that the resolvent of the NP operator is symmetric with respect to $0$. Inspired by computational results on rectangles we also prove that the NP operator on a rectangle has at least one eigenvalue if the aspect ratio of the rectangle is high enough. We first show that the spectral bound on rectangles tends to $1/2$ as the aspect ratio tends to $\infty$, and then existence of an eigenvalue follows as an immediate consequence. The numerical results of this paper also show that the interval determined by the bound on the essential spectrum obtained in \cite{PP-JAM-14} is actually the essential spectrum. After completion of the major part of this work we were informed by Mihai Putinar that he and Karl-Mikael Perfekt prove this rigorously \cite{PP-arXiv}. Their paper and the current one are complementary to each other.

The rest of the paper is organized as follows. In the next section we review symmetrization of the NP operator, and prove symmetry of spectrum and existence of eigenvalues on rectangles of high aspect ratio. In section \ref{sec:class} we present a method to classify spectra by quantifying resonance. For computation of resonance we use polarizable dipoles as a source function. We show advantages using polarizable dipoles in section \ref{sec:source}. Section \ref{sec:RCIP} is to describe the computational method (RCIP method) of this paper. Section \ref{sec:experi} is to present results of computational experiments. This paper ends with a short conclusion and discussion on mathematical problems raised by computational results.

\section{Spectrum of the NP operator in two dimensions}

Throughout this paper, we denote by $\la \cdot , \cdot \ra$ the duality  pairing of $H^{-1/2}$ and $H^{1/2}$, and $\| \cdot \|_{-1/2}$ denotes the $H^{-1/2}$ norm on $\p\GO$. Let $H_0^{-1/2}(\p{\GO})$ be the space of $\psi \in H^{-1/2}(\p{\GO})$ satisfying $\la \psi, 1 \ra =0$.

The single layer potential $\Scal_{\p\GO}[\Gvf]$ of a function $\Gvf$ on $\p\GO$ is defined by
\beq
\Scal_{\p\GO}[\Gvf](x) := \frac{1}{2\pi} \int_{\p\GO} \ln |x-y| \Gvf(y) \, d\Gs(y), \quad x \in \Rbb^2.
\eeq
Its relation to the NP operator is given by the following jump formula (see, for example, \cite{AK-book-07, Folland-book}):
\beq\label{j-single}
\p_\Gv \Scal_{\p\GO}[\Gvf] |_{\pm} = \left( \pm \frac{1}{2} I + \Kcal_{\p\GO}^* \right) [\Gvf] \quad \text{on } \p\GO,
\eeq
where $\p_\Gv$ denotes the outward normal derivative on $\p\GO$, and the subscripts $+$ and $-$ respectively indicate the limits (to $\p\GO$) from outside and inside $\GO$.

It is found in \cite{KPS-ARMA-07} that $\Kcal_{\p\GO}^*$ can be symmetrized using Plemelj's symmetrization principle
\beq\label{Plemelj}
\Scal_{\p\GO} \Kcal_{\p\GO}^* = \Kcal_{\p\GO} \Scal_{\p\GO}.
\eeq
If we define, for $\Gvf, \psi \in H_0^{-1/2}(\p{\GO})$,
\beq\label{newinner}
( \Gvf, \psi )_* := - \la \Gvf, \Scal_{\p\GO}[\psi] \ra = - \frac{1}{2\pi} \int_{\p\GO} \int_{\p\GO} \ln |x-y| \Gvf(x) \ol{\psi(y)} \, d\Gs(x)d\Gs(y),
\eeq
then $( \cdot, \cdot )_*$ is an inner product on $H_0^{-1/2}(\p{\GO})$, and the norm $\| \cdot \|_{*}$ induced by this inner product is equivalent to the $H^{-1/2}(\p\GO)$ norm, namely,
\beq
\| \Gvf \|_{\Hcal^*} \approx \| \Gvf \|_{-1/2}
\eeq
for all $\Gvf \in H_0^{-1/2}(\p{\GO})$ (see \cite{KKLS}).
Let $\Hcal^*_0$ be the space $H_0^{-1/2}(\p{\GO})$ equipped with the inner product $( \cdot, \cdot )_*$. Then the symmetrization principle \eqref{Plemelj} shows that $\Kcal_{\p\GO}^*$ is self-adjoint on $\Hcal^*_0$.

Let $\Gs(\Kcal_{\p\GO}^*)$ be the spectrum of $\Kcal_{\p\GO}^*$ on $\Hcal^*_0$. Since $\Kcal_{\p\GO}^*$ is self-adjoint on $\Hcal^*_0$, $\Gs(\Kcal_{\p\GO}^*)$ consists of continuous spectrum and pure point spectrum (eigenvalues), and continuous spectrum in turn consists of absolutely continuous spectrum and singularly continuous spectrum, namely,
\beq
\Gs(\Kcal_{\p\GO}^*)= \Gs_{\rm{c}}(\Kcal_{\p\GO}^*) \cup \ol{\Gs_{\rm{pp}}(\Kcal_{\p\GO}^*)} = \Gs_{\rm{ac}}(\Kcal_{\p\GO}^*) \cup \Gs_{\rm{sc}}(\Kcal_{\p\GO}^*) \cup \ol{\Gs_{\rm{pp}}(\Kcal_{\p\GO}^*)} ,
\eeq
and continuous spectrum and pure point spectrum are mutually disjoint (see \cite{RS-book-80, Yosida-book-74}). It is known (see \cite{Kellog-book}) that
\beq
\Gs(\Kcal_{\p\GO}^*) \subset (-1/2, 1/2).
\eeq

We will present a method based on resonance to distinguish continuous spectrum from pure point spectrum in the section \ref{sec:class}. Results of numerical experiments presented in section \ref{sec:experi} show that the spectrum is symmetric with respect to $0$. They also shows that on rectangles more and more eigenvalues appear as the aspect ratio increases, in particular, the pure point spectrum is non-empty if the aspect ratio is high enough. Let us prove rigorously these facts first in the following subsections.

\subsection{Symmetry of spectrum}

Here we prove that the spectrum of the NP operator in two dimensions is symmetric with respect to $0$, namely, that $\Gl$ is in the spectrum if and only if $-\Gl$ is. As mentioned before, this fact for eigenvalues is known \cite{Schiffer-PJM-57}. We extend it to include continuous spectrum. We emphasize that the spectrum we are considering is that on $\Hcal^*_0$. It is known that $1/2$ is an eigenvalue of $\Kcal_{\p\GO}^*$ on $H^{-1/2}(\p{\GO})$ while $-1/2$ is not. It is worth mentioning that the spectrum of the NP operator in three dimensions may not be symmetric with respect to $0$. For example, eigenvalues on the ball are all positive (see, for example, \cite{Kang-SMF-15, KPS-ARMA-07}).

We have the following theorem.

\begin{theorem}\label{thm:pair}
It holds in two dimensions that
\beq\label{cppsymm}
\Gs_{\rm{c}}(\Kcal_{\p\GO}^*) = - \Gs_{\rm{c}}(\Kcal_{\p\GO}^*), \quad \Gs_{\rm{pp}}(\Kcal_{\p\GO}^*) = - \Gs_{\rm{pp}}(\Kcal_{\p\GO}^*).
\eeq
\end{theorem}
\pf
For a given $\psi \in \Hcal_0^*$ let $u_\psi$ be the solution to
$$
\begin{cases}
\GD u_\psi=0 \quad &\text{in } \GO, \\
\p_\Gv u_\psi= \psi \quad &\text{on } \p\GO.
\end{cases}
$$
Let $u_\psi^\perp$ be a harmonic conjugate of $u_\psi$ in $\GO$ so that
\beq
\p_\Gj u_\psi^\perp= \p_\Gv u_\psi, \quad \p_\Gv u_\psi^\perp = - \p_\Gj u_\psi,
\eeq
where $\p_\Gj$ denotes the tangential derivative on $\p\GO$. Let $\psi^\perp:= \p_\Gv u_\psi^\perp = - \p_\Gj u_\psi$ on $\p\GO$, namely, $\psi^\perp$ is the Hilbert transform of $\psi$.

We first prove that
\beq\label{csymm}
\Gs(\Kcal_{\p\GO}^*) = - \Gs(\Kcal_{\p\GO}^*).
\eeq
To do so, it suffices to show that $\Gl$ is in the resolvent if and only if $-\Gl$ is. Suppose that $\Gl$ is in the resolvent, namely, $\Gl I - \Kcal_{\p\GO}^*$ is invertible on $\Hcal_0^*$.
Let $k$ be the number such that
\beq
\Gl= \frac{k+1}{2(k-1)}.
\eeq
For $\psi \in \Hcal_0^*$ let $\Gvf\in \Hcal_0^*$ be the unique solution to
\beq\label{intsol}
(\Gl I - \Kcal_{\p\GO}^*)[\Gvf]= \frac{1}{k-1} \psi^\perp.
\eeq
We see from \eqnref{j-single} that $u(x):= \Scal_{\p\GO}[\Gvf](x)$ for $x \in \Rbb^2$ is a solution to
\beq\label{usol}
\begin{cases}
\GD u= 0 \quad\text{in } \GO \cup (\Rbb^2 \setminus \ol{\GO}), \\
u|_+ - u|_- = 0, \\
\p_\Gv u|_+ - k \p_\Gv u|_- = \psi^\perp.
\end{cases}
\eeq
Moreover, since $\la \Gvf , 1 \ra=0$, it holds that
\beq\label{tendtozero}
u(x) \to 0 \quad \text{as } |x| \to \infty.
\eeq

Let $v_i$ be a harmonic conjugate of $u$ in $\GO$. Thanks to \eqnref{tendtozero}, $u$ has a harmonic conjugate in $\Rbb^2 \setminus \ol{\GO}$. Let $v_e$ be the harmonic conjugate in $\Rbb^2 \setminus \ol{\GO}$ such that $v_e(x) \to 0$ as $|x| \to \infty$. Then $v_i$ and $v_e$ are harmonic in $\GO$ and $\Rbb^2 \setminus \ol{\GO}$, respectively, and they satisfy
$$
\begin{cases}
\p_\Gj v_e - k \p_\Gj v_i = \p_\Gv u|_+ - k \p_\Gv u|_- = \psi^\perp \\
\p_\Gv v_e - \p_\Gv v_i = -\p_\Gj (u|_+ - u|_-) = 0
\end{cases}
\quad\text{on } \p\GO.
$$
Define $w$ by
\beq\label{wdef}
w(x):= \begin{cases}
k v_i(x)- u_\psi(x) +C, \quad & x \in \GO, \\
v_e(x), \quad & x \in \Rbb^2 \setminus \ol{\GO},
\end{cases}
\eeq
where $C$ is a constant to be determined. Then, $w$ is a solution to
\beq\label{wsol}
\begin{cases}
\GD w= 0 \quad\text{in } \GO \cup (\Rbb^2 \setminus \ol{\GO}), \\
w|_+ - w|_- = 0 \quad\text{on } \p\GO \\
\p_\Gv w|_+ - \frac{1}{k} \p_\Gv w|_- = \frac{1}{k} \psi \quad\text{on } \p\GO,  \\
w(x) \to 0 \quad\text{as } |x| \to \infty.
\end{cases}
\eeq
In fact, we have
$$
\p_\Gj w|_+ - \p_\Gj w|_- = \p_\Gj v_e - k \p_\Gj v_i + \p_\Gj u_\psi = \psi^\perp - \psi^\perp=0.
$$
So, $w|_+ - w|_-$ is constant on $\p\GO$. Hence, we can make it vanish on $\p\GO$ by choosing the constant $C$ properly.

Define $\Gvf_1 \in \Hcal_0^*$ by
\beq
\Gvf_1:= \p_\Gv w|_+ - \p_\Gv w|_-.
\eeq
Then we have
\beq\label{wSGvf}
w(x)= \Scal_{\p\GO}[\Gvf_1](x), \quad x \in \Rbb^2.
\eeq
In fact, if let $W:=w- \Scal_{\p\GO}[\Gvf_1]$, then we see from \eqnref{j-single} that
$$
\p_\Gv W|_+ - \p_\Gv W|_- = 0 \quad\text{on } \p\GO.
$$
So, it follows from \eqnref{wsol} that $W$ is a solution to
$$
\begin{cases}
\GD W= 0 \quad\text{in } \GO \cup (\Rbb^2 \setminus \ol{\GO}), \\
W|_+ - W|_- = 0 \quad\text{on } \p\GO \\
\p_\Gv W|_+ - \p_\Gv W|_- = 0 \quad\text{on } \p\GO,  \\
W(x) \to 0 \quad\text{as } |x| \to \infty.
\end{cases}
$$
It then follows from Green's identity that
$$
\int_{\GO} |\nabla W|^2 + \int_{\Rbb^2 \setminus \ol{\GO}} |\nabla W|^2 = \int_{\p\GO} W|_- (\p_\Gv W|_- - \p_\Gv W|_+ )=0,
$$
and hence $W$ is constant. Since $W(x) \to 0$ as $|x| \to \infty$, $W \equiv 0$. So we have \eqnref{wSGvf}.
By plugging \eqnref{wSGvf} into the third identity in \eqnref{wsol}, one can see from \eqnref{j-single} that
\beq\label{solvable}
(-\Gl I - \Kcal_{\p\GO}^*)[\Gvf_1]= \frac{1}{1-k} \psi.
\eeq

So far we have shown that $-\Gl I - \Kcal_{\p\GO}^*$ is surjective on $\Hcal_0^*$. Injectivity can be proved by reversing arguments from \eqnref{solvable} (with $\psi=0$ and $-\Gl$ replaced with $\Gl$) to \eqnref{intsol}. By interchanging the role of $\Gl$ and $-\Gl$ we see that $\Gl I - \Kcal_{\p\GO}^*$ is invertible if and only if $-\Gl I - \Kcal_{\p\GO}^*$ is. This proves \eqnref{csymm}.

Let us now prove \eqnref{cppsymm}. To do so, it suffices to prove the second identity
since $\Gs_{\rm{c}}(\Kcal_{\p\GO}^*) \cap \Gs_{\rm{pp}}(\Kcal_{\p\GO}^*) = \emptyset$. We include a short proof here even if it is a known fact as mentioned before.

Suppose that
\beq
(\Gl I - \Kcal_{\p\GO}^*)[\Gvf]= 0
\eeq
for some non-zero $\Gvf\in \Hcal_0^*$. Then, $u(x):= \Scal_{\p\GO}[\Gvf](x)$ is a nontrivial solution to
\eqnref{usol} with $\psi^\perp=0$ satisfying \eqnref{tendtozero}.
Then $w$ defined by \eqnref{wdef} is a solution to \eqnref{wsol} with $\psi=0$. Then, $w(x)= \Scal_{\p\GO}[\Gvf_2](x)$ for some non-zero
$\Gvf_2 \in \Hcal_0^*$, and it holds that
\beq
(-\Gl I - \Kcal_{\p\GO}^*)[\Gvf_2]= 0.
\eeq
So we have shown that if $\Gl \in \Gs_{\rm{pp}}(\Kcal_{\p\GO}^*)$, then $-\Gl \in \Gs_{\rm{pp}}(\Kcal_{\p\GO}^*)$. By interchanging the role of $\Gl$ and $-\Gl$ we show the second identity in \eqnref{cppsymm}.
This completes the proof. \qed

\subsection{Existence of eigenvalues on rectangles of high aspect ratio}

Let us first recall that the spectral bound $b$ (other than $1/2$) of the NP operator on $\p\GO$ is given by
\beq\label{specbound}
b = \frac{1}{2} \sup_{\Gvf \in \Hcal^*_0} \frac{\left| \| \nabla \Scal_{\p\GO}[\Gvf] \|_{L^2(\Rbb^d \setminus \GO)}^2 -
\| \nabla \Scal_{\p\GO}[\Gvf] \|_{L^2(\GO)}^2 \right|}{\| \nabla \Scal_{\p\GO}[\Gvf] \|_{L^2(\Rbb^d)}^2}.
\eeq
(See, for example, \cite{KPS-ARMA-07}.) On the other hand, the bound on the essential spectrum of the NP operator on curvilinear polygonal domains is obtained in \cite{PP-JAM-14}:
\beq\label{bounds}
b_{\text{ess}}= \frac{1}{2} \max_{1\le j \le N} \left(1- \frac{\Gt_j}{\pi} \right),
\eeq
where $\Gt_j$ is the interior angle of the $j$th corner and $N$ is the number of corners. Note that $1/2$ appears in \eqnref{bounds} since the NP operator of this paper is $1/2$ times the one in \cite{PP-JAM-14}. If $\GO$ is a rectangle, $b_{\text{ess}}=1/4$. We show that on rectangles of high aspect ratio $b$ is larger than $1/4$, and hence there must be eigenvalues.

\begin{theorem}\label{thm:limbr}
For $r \ge 1$, let $\GO_r$ be a rectangle of aspect ratio $r$ and $b_r$ be the spectral bound of the NP operator on $\GO_r$. It holds that
\beq\label{limbr}
\lim_{r \to \infty} b_r = \frac{1}{2}.
\eeq
\end{theorem}

Above theorem shows that there is $r_0$ such that if $r\ge r_0$, then $b_r>1/4$. It means that there is a member of spectrum $\Gl$ such that $1/4< \Gl$. Since the essential spectrum is confined in $[-1/4, 1/4]$ due to \eqnref{bounds}, $\Gl$ must be an eigenvalue. So we have the following corollary.

\begin{cor}\label{cor:thinrec}
There is $r_0$ such that for any $r \ge r_0$ the NP operator on $\GO_r$ has at least one eigenvalue.
\end{cor}

\noindent{\sl Proof of Theorem \ref{thm:limbr}}.
We adapt the idea of the proof of Theorem 5 in \cite{KPS-ARMA-07}. For $\Gd >0$ there exists a $C^\infty$ odd function $\psi$ on $\Rbb$ such that \begin{align}
& \psi(t) \le C \Gd^{-1} \quad\text{for all } t, \label{psibound} \\
& \psi(t)=t \quad\text{if } -1 \le t \le 1, \label{psit} \\
& \int_{|t| \ge 1} |\psi'(t)|^2 dt < \Gd, \label{psiint}
\end{align}
where $C$ is a universal constant independent of $\Gd$. Existence of such a function is proved in \cite[Lemma 7]{KPS-ARMA-07}. In fact, the property \eqnref{psibound} is not presented there. But one can easily check that $\psi$ there satisfies \eqnref{psibound}. Let $\chi$ be a non-negative $C^\infty$ function on $\Rbb$ with a compact support such that $\chi(t)=1$ on $[-1,1]$. Let $\chi_N(t):= \chi(t/N)$. Then, there is a constant $C$ independent of $N$ such that
\beq\label{chibound}
\chi_N(t) \le C \quad\text{for all } t,
\eeq
and
\beq\label{chiint}
\int_\Rbb \chi_N(t)^2 dt \le CN, \quad  \int_\Rbb \chi_N'(t)^2 dt \le CN^{-1}.
\eeq

Let us use $(x,y)$ for Cartesian coordinates in this proof. Since the NP operator is scale invariant, we may assume $\GO=\GO_r:= [-1,1] \times [-\Ge, \Ge]$ with $\Ge=r^{-1}$. Let $\psi_N(t):= \chi_N(t) \psi(t)$ and define
$$
w(x,y):= \chi_N(x) \psi_M(y/\Ge),
$$
where $M$ and $N$ are large numbers to be determined. The function $w$ also depends on $\Gd$ which is also to be determined. Note that $w(x,y)= y/\Ge$ if $(x,y) \in \GO$, in particular, it is harmonic there. Define
$$
\Gvf:= (-\frac{1}{2} I + \Kcal_{\p\GO}^*)^{-1} [\p_\nu w|_{\p\GO}] \quad \text{on } \p\GO,
$$
and let
$$
u(x,y):= \Scal_{\p\GO} [\Gvf](x,y), \quad (x,y) \in \Rbb^2.
$$
We emphasize that $\Gvf$ belongs to $\Hcal_0^*$. Since $\p_\nu u|_{-}= \p_\nu w$ on $\p\GO$, we have $u=w+C$ in $\GO$ for some constant $C$. One can easily see that $u$ is odd with respect to $x$-axis, so $C=0$. It then follows that
\beq\label{estGO}
\| \nabla u \|_{L^2(\GO)}^2 = \| \nabla w \|_{L^2(\GO)}^2 = 4 \Ge^{-1}.
\eeq

Let us now estimate $\| \nabla u \|_{L^2(\Rbb \setminus \GO)}^2$. Since $u=w$ on $\p\GO$, we have from Dirichlet's principle
$$
\| \nabla u \|_{L^2(\Rbb \setminus \GO)}^2 \le \| \nabla w \|_{L^2(\Rbb \setminus \GO)}^2= \int_{-\infty}^{\infty} \int_{|x| \ge 1} + \int_{|y| \ge \Ge} \int_{|x| \le 1} |\nabla w|^2 dxdy.
$$
So, we have
\begin{align*}
\| \nabla u \|_{L^2(\Rbb \setminus \GO)}^2 &\le \int_{-\infty}^{\infty} \left| \psi_M \left(\frac{y}{\Ge} \right) \right|^2 dy
\int_{|x| \ge 1} |\chi_N'(x)|^2 dx + \frac{1}{\Ge^2} \int_{|y| \ge \Ge} \left| \psi_M' \left(\frac{y}{\Ge} \right) \right|^2 dy \int_{|x| \ge 1} |\chi_N(x)|^2 dx \\
& \quad + \frac{1}{\Ge^2} \int_{|y| \ge \Ge} \left| \psi_M' \left(\frac{y}{\Ge} \right) \right|^2 dy \int_{|x| \le 1} |\chi_N(x)|^2 dx \\
&= \Ge \int_{-\infty}^{\infty} | \psi_M (y)|^2 dy
\int_{|x| \ge 1} |\chi_N'(x)|^2 dx + \frac{1}{\Ge} \int_{|y| \ge 1} |\psi_M' (y)|^2 dy \int_{|x| \ge 1} |\chi_N(x)|^2 dx \\
& \quad + \frac{1}{\Ge} \int_{|y| \ge 1} | \psi_M' (y) |^2 dy \int_{|x| \le 1} |\chi_N(x)|^2 dx
=: I_1+I_2+I_3.
\end{align*}
One can see from \eqnref{psibound} and the first inequality in \eqnref{chiint} that
$$
\int_{-\infty}^{\infty} | \psi_M (y)|^2 dy \le C \Gd^{-1} \int_{-\infty}^{\infty} | \chi_M (y)|^2 dy \le C \Gd^{-1} M.
$$
So we obtain from the second inequality in \eqnref{chiint} that
\beq\label{intone}
I_1 \le C \Ge \Gd^{-1} M N^{-1}.
\eeq
Similarly, we have
$$
\int_{|y| \ge 1} |\psi_M' (y)|^2 dy \le C \left( \Gd^{-1} \int_{|y| \ge 1} |\chi_M' (y)|^2 dy + \int_{|y| \ge 1} |\psi' (y)|^2 dy \right) \le C(\Gd^{-1} M^{-1} +\Gd),
$$
and hence
\beq\label{inttwo}
I_2 \le C \Ge^{-1} (\Gd^{-1} M^{-1} +\Gd) N,
\eeq
and
\beq\label{intthree}
I_3 \le C \Ge^{-1} (\Gd^{-1} M^{-1} +\Gd) .
\eeq
Putting them together, we see that
$$
\| \nabla u \|_{L^2(\Rbb \setminus \GO)}^2 \le C \left( \Ge \Gd^{-1} M N^{-1} + \Ge^{-1} \Gd^{-1} M^{-1} N + \Ge^{-1} \Gd N \right).
$$

Let us take, for example, $\Gd= \Ge^{1/3}$, $M=\Ge^{-2/3}$, and $N=\Ge^{-1/6}$. Then we have
$$
\| \nabla u \|_{L^2(\Rbb \setminus \GO)}^2 \le C \Ge^{-5/6}.
$$
According to \eqnref{specbound} and \eqnref{estGO}, the spectral bound $b_r$ satisfies
$$
b_r \ge \frac{1}{2} \frac{4 \Ge^{-1} - C \Ge^{-5/6}}{4 \Ge^{-1} + C \Ge^{-5/6}}.
$$
Since $b_r < 1/2$, we have \eqnref{limbr}. This completes the proof.
\qed

\section{Classification of spectrum by resonance}\label{sec:class}

Let $f \in \Hcal_0^*$. For $t \in (-1/2, 1/2)$ and $\Gd >0$, let $\Gvf_{t, \Gd}$ be the solution of the integral equation
\beq\label{inteqn}
\big( (t+i\Gd)I - \Kcal_{\p\GO}^* \big)[\Gvf_{t,\Gd}]= f \quad\mbox{on } \p\GO.
\eeq
By spectral resolution theorem \cite{Yosida-book-74}, there is a family of projection operators
$\Ecal_s$ (called the resolution identity) such that
\beq
\Kcal^*_{\p\GO} = \int_{\Gs(\Kcal_{\p\GO}^*)} s \, d  \Ecal_s.
\eeq
We then obtain from \eqnref{inteqn}
\beq\label{specGvf}
\Gvf_{t,\Gd} = \int_{\Gs(\Kcal_{\p\GO}^*)} \frac{1}{t+i\Gd -s} \, d \Ecal_s[f],
\eeq
and hence
\beq\label{GvfGd}
\| \Gvf_{t,\Gd} \|_{*}^2 =\int_{\Gs(\Kcal_{\p\GO}^*)} \frac{1}{(s-t)^2 + \Gd^2} \,d ( f, \Ecal_s[f] )_*.
\eeq

If $t \notin \Gs(\Kcal_{\p\GO}^*)$, one can immediately see from \eqnref{GvfGd} that
\beq\label{notspec}
\lim_{\Gd \to 0} \Gd^\Ga \| \Gvf_{t, \Gd} \|_{*} =0
\eeq
for any $\Ga>0$.

Suppose that $t \in \Gs(\Kcal_{\p\GO}^*)$. An eigenvalue $t$ of $\Kcal_{\p\GO}^*$ is characterized by discontinuity $\Ecal_{t+} - \Ecal_t \neq 0$ (and $t$ is isolated). So, if $f$ satisfies
\beq
( f, \Ecal_{t+}[f] )_* - ( f, \Ecal_t[f] )_* >0,
\eeq
then
$$
\| \Gvf_{t,\Gd} \|_{*}^2  \ge  \frac{( f, \Ecal_{t+}[f] )_* - ( f, \Ecal_t[f] )_* }{\Gd^2},
$$
and hence
\beq\label{eigenvalue}
\| \Gvf_{t,\Gd} \|_{^*}^2 \approx \Gd^{-2}.
\eeq

Suppose that the spectral measure $d ( f, \Ecal_s[f] )_*$ is absolutely continuous near $t$, namely, there is $\Ge >0$ and a function $\Gm_f(s)$ which is integrable on $[t-\Ge, t+\Ge]$ such that
\beq
d ( f, \Ecal_s[f] )_* = \Gm_f(s) ds, \quad s \in [t-\Ge, t+\Ge].
\eeq
Then we obtain from \eqnref{GvfGd}
$$
\| \Gvf_{t,\Gd} \|_{*}^2 =\int_{\Gs(\Kcal_{\p\GO}^*) \setminus [t-\Ge, t+\Ge]} \frac{1}{(s-t)^2 + \Gd^2} \,d ( f, \Ecal_s[f] )_* + \int_{[t-\Ge, t+\Ge]} \frac{\Gm_f(s) ds}{(s-t)^2 + \Gd^2} .
$$
From the boundary behavior of the Poisson integral, we have
\beq
\lim_{\Gd \to 0} \Gd \int_{[t-\Ge, t+\Ge]} \frac{\Gm_f(s) ds}{(s-t)^2 + \Gd^2} = \frac{\Gp}{2} (\Gm_f(t+)+\Gm_f(t-)).
\eeq
On the other hand, it is proved in \cite{KLY} that
\beq
\lim_{\Gd \to 0} \Gd^2 \int_{[t-\Ge, t+\Ge]} \frac{\Gm_f(s) ds}{(s-t)^2 + \Gd^2} = 0.
\eeq
So, we have
\beq\label{eqn:deltabeta}
\lim_{\Gd \to 0} \Gd \| \Gvf_{t,\Gd} \|_{*}^2 = \frac{\Gp}{2} (\Gm_f(t+)+\Gm_f(t-)),
\eeq
and
\beq\label{onezero}
\lim_{\Gd \to 0} \Gd^2 \| \Gvf_{t,\Gd} \|_{*}^2 = 0.
\eeq

Define an indicator function $\Ga_f(t)$ by
\beq
\Ga_f(t):= \sup \left\{ ~ \Ga ~\Big| ~ \limsup_{\Gd \to 0} \Gd^\Ga \| \Gvf_{t, \Gd} \|_{*} = \infty ~ \right\}, \quad t \in (-1/2, 1/2).
\eeq
We see that $0 \le \Ga_f(t) \le 1$ for all $t$. The classification of spectra of the NP operator is based on the following theorem.

\begin{theorem}\label{thm:class}
Let $f\in \Hcal_0^*$.
\begin{itemize}
\item[{\rm (i)}] If $\Ga_f(t) >0$, then $t \in \Gs(\Kcal_{\p\GO}^*)$.
\item[{\rm (ii)}] If $\Ga_f(t) =1$ and $t$ is isolated, then $t \in \Gs_{\rm{pp}}(\Kcal_{\p\GO}^*)$.
\item[{\rm (iii)}] If $1/2 \le \Ga_f(t) <1 $, then $t \in \Gs_{\rm{c}}(\Kcal_{\p\GO}^*)$.
\end{itemize}
\end{theorem}
\pf
The assertion (i) is an immediate consequence of \eqnref{notspec}, while (ii) follows from \eqnref{eigenvalue} and \eqnref{onezero}. (iii) is a consequence of (i) and (ii).
\qed

Because of complicated nature of the singularly continuous spectrum, it is hard to classify it from the continuous spectrum. However, $\Ga_f(t) =1$ and $t$ is not isolated (embedded in the continuous spectrum), then we may infer that it is in the singularly continuous spectrum.

The indicator function $\Ga_f(t)$ can be computed using the following lemma.

\begin{lemma}\label{lem:alpha}
For $t \in (-1/2, 1/2)$ and $\Gd >0$ define
\beq
\Ga_f(t,\Gd)=- \frac{\log \| \Gvf_{t, \Gd} \|_{*}}{\log \Gd}.
\eeq
It holds that
\beq
\Ga_f(t) =  \lim_{\Gd \to 0}\Ga_f(t,\Gd)
\eeq
if the limit exists.
\end{lemma}

\pf
Let $\ell:= \lim_{\Gd \to 0}\Ga_f(t,\Gd)$ under the assumption that the limit exists. If $\Ga < \ell$, then
$$
\limsup_{\Gd \to 0} \Gd^\Ga \| \Gvf_{t, \Gd} \|_{*} = \infty,
$$
and hence $\Ga \le \Ga_f(t)$. If $\Ga > \ell$, then
$$
\lim_{\Gd \to 0} \Gd^\Ga \| \Gvf_{t, \Gd} \|_{*} = 0,
$$
and hence $\Ga \ge \Ga_f(t)$. So $\Ga_f(t)=\ell$. \qed

Lemma \ref{lem:alpha} allows us to approximate $\Ga_f(t)$ by $\Ga_f(t,\Gd)$ for small $\Gd$. In fact, the high precision method to compute $\Kcal_{\p\GO}^*$, which will be described in the next section, makes it possible to use $\Gd$ smaller than $10^{-10}$.

\section{Source functions}\label{sec:source}

It is quite important for the classification of spectra using $\Ga_f$ to choose the source function properly. For example, to have $\Ga_f(t)=1$ for classification of eigenvalues, $f$ needs to have a non-zero eigenmode (the component of the corresponding eigenfunctions). To be more precise, we consider \eqnref{inteqn} when $\Gd=0$, namely,
\beq\label{inteqn2}
\big( t I - \Kcal_{\p\GO}^* \big)[\Gvf]= f \quad\mbox{on } \p\GO .
\eeq
Since $\Kcal_{\p\GO}^*$ is self-adjoint, this equation is solvable if and only if $f \perp \text{Ker} (t I - \Kcal_{\p\GO}^*)$. It means that in order to characterize spectrum of $\Kcal_{\p\GO}^*$ in terms of blow-up of $\| \Gvf_{t, \Gd} \|_*$, $f$ must have a component of $\text{Ker} (t I - \Kcal_{\p\GO}^*)$.

For example, if $f=\nu \cdot \nabla (d \cdot x)$ for a constant vector
$d$, there are dark plasmons which are the eigenvalues undetectable by
$f$ \cite{HP-ACHA-13}. In this regard, it is helpful to mention about
the polarization tensors. Let $\Gvf_{t, \Gd}^{(j)}$ be the solution
\eqnref{inteqn} when $f=\nu \cdot \nabla x_j$. For $i,j=1,2$, we
define
$$
m_{ij} (t+i\Gd):= \int_{\p\GO} x_i \Gvf_{t, \Gd}^{(j)}(x) \, d\Gs.
$$
The matrix $M(t+i\Gd):= (m_{ij} (t+i\Gd))$ is called the polarization tensor. It is an analytic function of $\Gl=t+i\Gd$ in $\Cbb \setminus (-1/2, 1/2)$, and may have singularities at $\Gl=t$ in the spectra of $\Kcal_{\p\GO}^*$. The singularities of $M(t+i\Gd)$ were investigated in \cite{HMM-NJP-11, HP-ACHA-13} when $M$ is isotropic. However, singularities of $M(t+i\Gd)$ can show some spectrum, but not all as the following example shows. If $\GO$ is an ellipse of major axis $a$ and minor axis $b$, then $M(\Gl)=M(t+i\Gd)$ is given by
\beq\label{PTellipse}
M(\Gl) = 2 \pi ab
\begin{bmatrix}
 \ds  \frac{a+b}{(2\Gl-1)a+ (2\Gl+1)b} & 0 \\ \nm
  0 & \ds \frac{a+b}{(2\Gl-1)b + (2\Gl+1)a}
\end{bmatrix}.
\eeq
(See, for example, \cite{Kang-SMF-15}.) So the singularities (actually poles) of $M(\Gl)$ occur only at $\pm \frac{a-b}{2(a+b)}$. However, it is known that eigenvalues of $\Kcal_{\p\GO}^*$ are
\beq\label{eqn:ellipse}
\pm \frac{1}{2} \left( \frac{a-b}{a+b} \right)^n, \quad n=1,2,\ldots.
\eeq
So, in this case $M(\Gl)$ shows only two eigenvalues. It is quite interesting to observe that those two eigenvalues are the largest (in absolute values) eigenvalues.

In this paper we use as source functions \beq\label{eqn:f} f_z(x)
=\nu(x)\cdot\nabla q_z(x),\quad q(x)=d\cdot \nabla_x \left(
  \frac{1}{2\pi}\ln|x-z|\right) \eeq where $d$ is a constant unit
vector and $z\in \Rbb^2\setminus\overline{\Omega}$. In fact, $q_z$ is
the newtonian potential of the polarizable dipole
$d\cdot\nabla\Gd_z(x)$ ($\Gd_z (x)$ is the Dirac mass) located at $z$.
The source function $f_z$ was used in \cite{KLY} for analysis of
resonance on intersecting disks.

Using $f_z$ as source functions has several advantages. First of all, $q_z$ is a harmonic function in $\GO$, and hence
$$
\| f_z \|_{H^{-1/2}(\p\GO)} \approx \| q_z \|_{H^1(\GO)}
$$
by the standard regularity estimates of the Neumann problem for the Laplace equation. Moreover, one can see easily that there are constants $C_1$ and $C_2$ such that
$$
C_1 \le \| q_z \|_{H^1(\GO)} \le C_2
$$
as long as the location $z$ of the dipole is at some distance from $\GO$ and $\infty$, namely, there are constants $C_3$ and $C_4$ such that
\beq\label{distance}
C_3 \le \text{dist}(z, \GO) \le C_4.
\eeq
It means that for all $z$ satisfying \eqnref{distance} we have
\beq
\| f_z \|_{H^{-1/2}(\p\GO)} \approx 1.
\eeq
In particular, we don't have to normalize $f_z$.

Another advantage of using $f_z$ as source functions is that for any $\Gvf \in \Hcal_0^*$ $f_z$ contains a component of $\Gvf$, namely, $(f_z, \Gvf)_* \neq 0$ for most $z$'s. To see this we first observe that
\begin{align*}
(f_z, \Gvf)_* &= - \la f_z, \Scal_{\p\GO}[\Gvf] \ra
= - \int_{\GO} \nabla q_z(x) \cdot \nabla \Scal_{\p\GO}[\Gvf](x) \, dx \\
&= d \cdot \nabla_z \left( \frac{1}{2\pi} \int_{\GO} \frac{x-z}{|x-z|^2} \cdot \nabla \Scal_{\p\GO}[\Gvf](x) \, dx \right) = d \cdot \nabla_z \Scal_{\p\GO} [\p_\nu \Scal_{\p\GO}[\Gvf]](z).
\end{align*}
It shows that $(f_z, \Gvf)_*$ is harmonic as a function of the $z$ variable, and it is non-vanishing. In fact, if $(f_z, \Gvf)_*=0$ for all $z \in \Rbb^2 \setminus \ol{\GO}$, then $\Scal_{\p\GO} [\p_\nu \Scal_{\p\GO}[\Gvf]](z)=0$, and hence $\p_\nu \Scal_{\p\GO}[\Gvf]=0$ on $\p\GO$. Thus $\Gvf=0$. As a non-vanishing harmonic function, $(f_z, \Gvf)_*$ cannot be zero for $z$ in an open set. So, we infer that $(f_z, \Gvf)_* \neq 0$ for almost all $z$.

Still $(f_z, \Gvf)_*$ can be small for $z$ in a large set. To avoid such a case, we choose several $z$'s, say $z_1, \ldots, z_N$, satisfying \eqnref{distance}, and consider the new indicator function
\beq\label{eqn:alphasharp1}
\Ga_\sharp(t):=\max_{1\le j \le N} \{ \Ga_{f_{z_j}}(t) \}.
\eeq
We emphasize that Theorem \ref{thm:class} is still valid with $\Ga_\sharp(t)$ replaced with $\Ga_f(t)$.
The indicator function $\Ga_\sharp(t)$ is approximated by $\Ga_\sharp(t,\Gd)$ for small $\Gd$, which is defined by
\beq\label{eqn:alphasharp2}
\Ga_\sharp(t,\Gd):= \max_{1\le j \le N} \{ \Ga_{f_{z_j}}(t, \Gd) \}.
\eeq

\section{Description of the numerical method}\label{sec:RCIP}

This section briefly motivates and discusses the numerical method used
in section 6 to solve \eqnref{inteqn} and to compute the inner
product $( \cdot, \cdot )_*$ in \eqnref {newinner}.

\subsection{Nystr\"{o}m discretization on smooth boundaries}\label{sec:Nystr}

Let $K(x,y)$ be the kernel of an integral operator $K$ that is compact
on a smooth boundary $\partial\Omega$ and let $f(x)$ be a smooth
function on $\partial\Omega$. A popular method for finding numerical
solutions to Fredholm second kind boundary integral equations of the
type
\begin{equation}
\varphi(x)+\int_{\partial\Omega} K(x,y)\varphi(y)\,d\Gs(y)=f(x)\,,
\quad x\in\partial\Omega\,,
\label{eq:inteq1}
\end{equation}
is {\it Nystr\"{o}m discretization}: the integral in~(\ref{eq:inteq1}) is
discretized on a mesh on $\partial\Omega$ according to some
polynomial-based quadrature rule with a number $N$ of nodes and
weights $x_j$ and $w_j$, $j=1,\ldots,N$, and the resulting
semi-discrete equation for the unknown layer density $\varphi(x)$ is
enforced at the quadrature nodes~\cite[Chapter~4.1]{Atki97}. Upon
solving the resulting linear system
\begin{equation}
\tilde{\varphi}(x_j)+
\sum_{k=1}^N K(x_j,x_k)\tilde{\varphi}(x_k)w_k=r(x_j)\,,\quad j=1,\ldots,N\,,
\label{eq:inteq2}
\end{equation}
one obtains an approximation $\tilde{\varphi}(x_j)$ to $\varphi(x_j)$ whose
convergence with $N$ reflects that of the underlying quadrature.
When~(\ref{eq:inteq1}) stems from a well-conditioned boundary value
problem on a domain whose boundary $\partial\Omega$ can be resolved
with a moderate number of discretization points, uniform meshes and
high-order accurate quadrature rules, such as composite 16-point
Gauss--Legendre quadrature, are appropriate in the sense that they
often produce solutions $\tilde{\varphi}(x_j)$ with a relative accuracy
close to machine epsilon ($\epsilon_{\rm mach})$ at modest
computational costs.

\subsection{Difficulties related to piecewise smooth boundaries}\label{subsection:num_diffi}

Eq.\;\eqnref{inteqn} of the present paper is not quite of the
type~(\ref{eq:inteq1}). The chief difference being that the boundary
$\partial\Omega$ has a finite number of corners where the NP operator
${\cal K}^*_{\partial\Omega}$ in \eqnref{inteqn} is not compact. This
lack of compactness manifests itself in that the solution
$\varphi_{t,\delta}(x)$ to \eqnref{inteqn} may exhibit a non-smooth,
oscillatory, and diverging behavior close to the corner vertices which
cannot easily be resolved by polynomials on a uniform mesh. Rather,
intense mesh refinement is needed for accuracy. This, in turn, may
lead to all kinds of numerical problems related to computational
economy and to stability irrespective of what numerical method is
used.

Another difficulty with producing numerical solutions
$\varphi_{t,\delta}(x)$ to \eqnref{inteqn} appears for very small
values of $\delta$ in combination with $t$ being close to, or in, the
spectrum of the NP operator. Finite precision arithmetic makes it hard
for any solver to distinguish between combinations of $\delta$ and $t$
for which a solution exists and combinations for which no solution
exists. Such dichotomies often imply numerical ill-conditioning and
the loss of precision.

Finally, the computation of $( \cdot, \cdot )_*$ offers challenges
when $\varphi_{t,\delta}(x)$ is non-smooth, oscillatory, and diverging.
The action of the operator ${\cal S}_{\partial\Omega}$ on
$\varphi_{t,\delta}(x)$ may become inaccurate due to numerical
cancellation even if $\varphi_{t,\delta}(x)$ itself is accurate.

\subsection{RCIP acceleration and Nystr\"{o}m schemes}\label{sec:RaN}

Fortunately, most of the numerical difficulties discussed in section
\ref{subsection:num_diffi}, can be overcome by the use of the recursively compressed inverse preconditioning (RCIP)
method~\cite{Helsing-tutorial,HelsOjal08}, which is a tool to improve
the stability and greatly reduce the computational cost of standard
Nystr\"{o}m discretization schemes when applied to Fredholm second
kind integral equations on piecewise smooth domains. In particular,
for $t$ at some distance away from zero, from the endpoints of the
continuous spectrum, and from the pure point spectrum,
RCIP-accelerated Nystr\"{o}m solvers in combination with fixed-point
iteration, Newton's method, and a certain homotopy
technique~\cite[Section 6]{Hels11} can often produce solutions
$\varphi_{t,\delta}(x)$ to \eqnref{inteqn} with a relative precision
of about $10\cdot\epsilon_{\rm mach}$ for any $\delta$. See
\cite[Fig.\;10]{Hels11} for an illustration. The achievable relative
precision for the inner product $( \cdot, \cdot )_*$ is, typically,
$10\cdot\epsilon_{\rm mach}/\delta$.

In the numerical examples of section \ref{sec:experi} in this paper we
use an RCIP-accelerated Nystr\"{o}m solver that has previously been
used to compute, very accurately, polarizabilities of various
arrangements of dielectric squares and
cubes~\cite{Helsing-tutorial,HP-ACHA-13} as well as electromagnetic
resonances inside microwave cavities with sharp
edges~\cite{HelsKarl16}. It would carry too far to recapitulate and
put in context the fairly large collection of numerical techniques
that constitute the RCIP method, so we refer the reader to
\cite{Hels11,Helsing-tutorial} for details.

\section{Numerical experiments}\label{sec:experi}

This section presents numerical results for the spectrum of the NP
operator on various domains as revealed by the indicator function
$\Ga_\sharp(t)$ of \eqnref{eqn:alphasharp1}. On smooth domains, the
spectrum consists only of pure point spectrum and the standard
Nystr{\"o}m method of section~\ref{sec:Nystr} is very efficient. This
is illustrated with experiments on ellipses and superellipses in
subsections~\ref{subsec:ell} and~\ref{subsec:rec}.  All other
experiments apply to domains with corners and use the RCIP-accelerated
Nystr{\"o}m solver mentioned in section~\ref{sec:RaN}. The high
performance of this solver is demonstrated by comparison with
analytical results for intersecting disks, which is the only domain
with corners for which the spectrum is fully analyzed~\cite{KLY}. We
present spectrum on a triangle, on rectangles of various aspect
ratios, and on an ellipse perturbed by a corner. The results for
rectangles show that the pure point spectrum can be empty or non-empty
depending on the aspect ratio.  The result for the perturbed ellipse shows
that there are eigenvalues embedded in the absolutely continuous
spectrum. So, we conclude that pure point spectrum and singularly
continuous spectrum can be non-empty depending on the geometry of the
domain.

From now on, let $f_z$ be the dipole source function in
\eqnref{eqn:f}. To determine $\Ga_\sharp(t)$ of
\eqnref{eqn:alphasharp1} for a particular domain we use $N=10^3$
source functions $f_{z_j}$ located on a circle enclosing the domain.
In addition to maximizing over dipole source locations we also, for
each $f_{z_j}$, maximize over $10^3$ orientations of the unit vector
$d$. In this way, the maximum in \eqnref{eqn:alphasharp1} is taken
over $10^6$ dipole fields for each $t$-value. The enclosing circle is
centered at the origin and has radius $R$.

Our experiments confirm the symmetry of Theorem~\ref{thm:pair}. They
also show that the bounds on the essential spectrum obtained
in~\cite{PP-JAM-14} are optimal. Actually, the experiments show even
more: the whole interval between the bounds is the essential spectrum. As mentioned in Introduction, this was proved rigorously quite lately in
\cite{PP-arXiv}.

\subsection{Ellipses}\label{subsec:ell}

\begin{figure}[!ht]
\begin{center}
\epsfig{figure=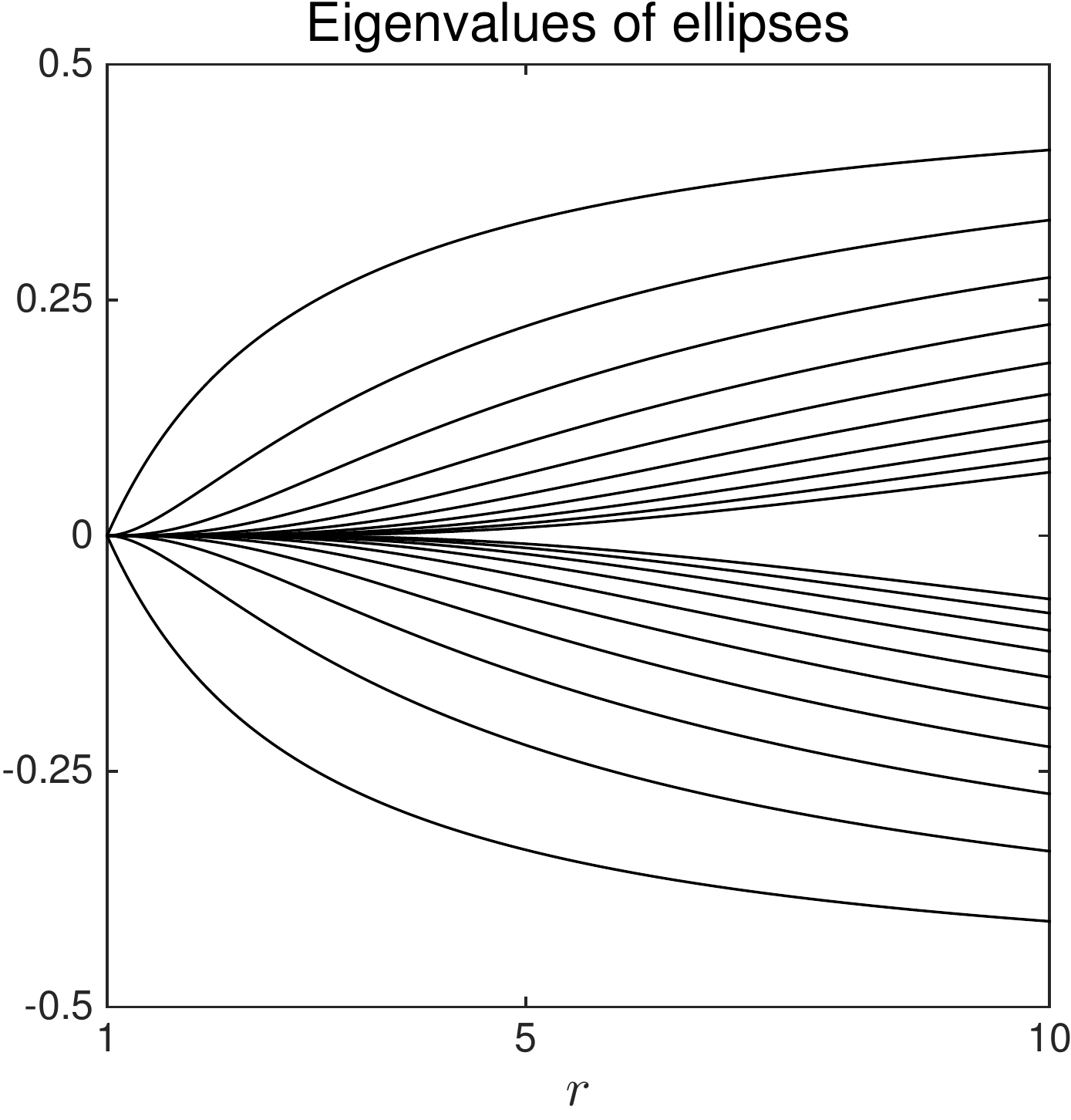,width=6cm}
\end{center}
\label{intersecting}
\caption{The 20 largest eigenvalues for ellipses of various aspect
  ratios $r$. Thinner ellipses have bigger eigenvalues.
  }\label{fig:ellipse}
\end{figure}

Recall from \eqnref{eqn:ellipse} that for an ellipse with aspect ratio
$r$, the eigenvalues of the NP operator are
\begin{displaymath}
\pm \frac{1}{2}\left(\frac{1-r}{1+r}\right)^n,\quad n=1,2,\dots.
\end{displaymath}
For $n\leq 40$ and $r=3$, the Nystr{\"o}m method reproduces these
eigenvalues with an absolute error of less than $7\cdot10^{-16}$. For
$r=30$ the error is less than $6\cdot10^{-15}$.
Fig.\;\ref{fig:ellipse} illustrates the largest eigenvalues for $1\leq
r\leq 10$.

\subsection{Intersecting disks}\label{subsec:inter}

\begin{figure}[htp]
    \begin{subfigure}{0.48\textwidth}
      \centering
      \includegraphics[height=5cm]{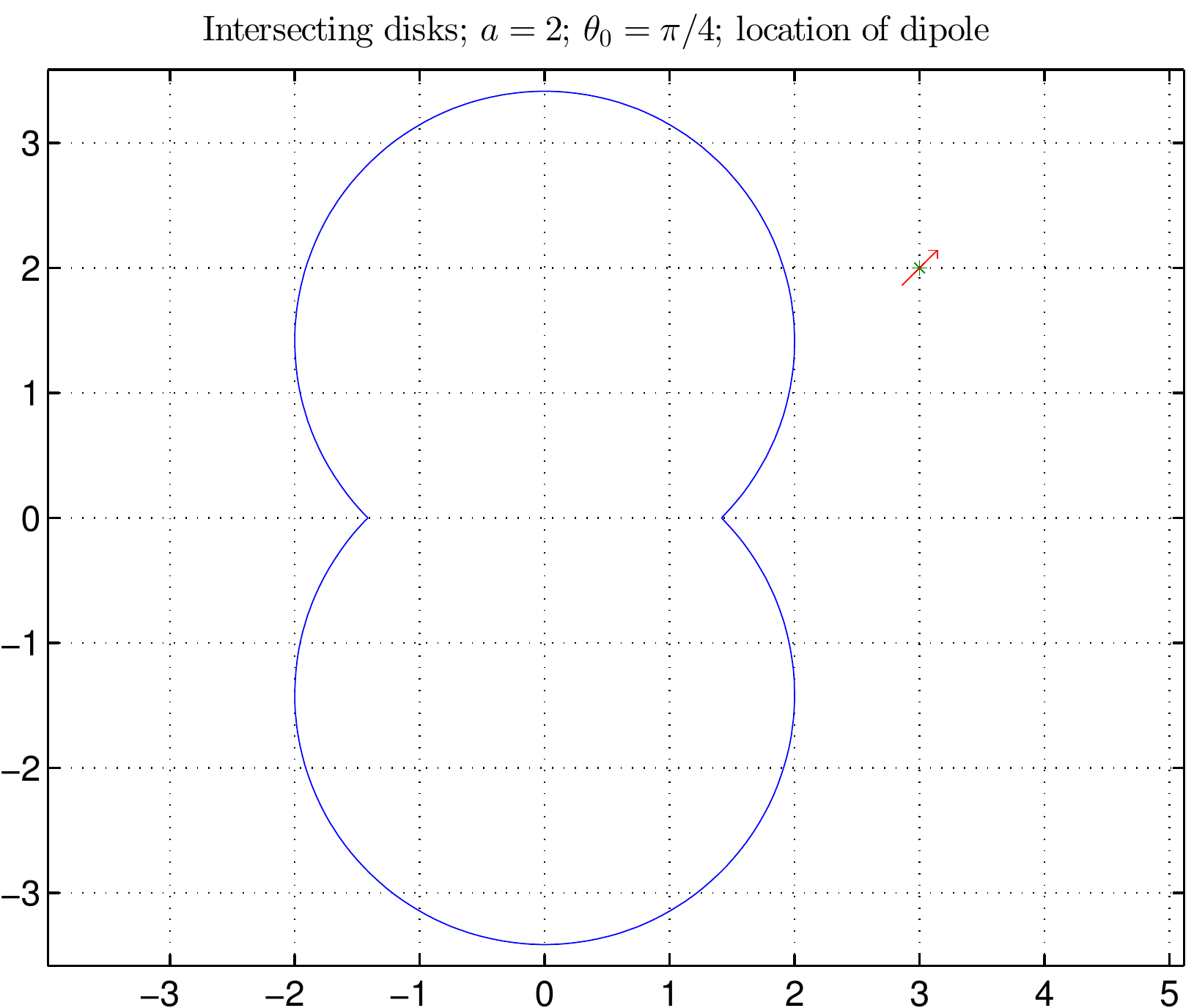}
      \caption{}
    \end{subfigure}
       \begin{subfigure}{0.48\textwidth}
      \centering
      \includegraphics[height=5cm]{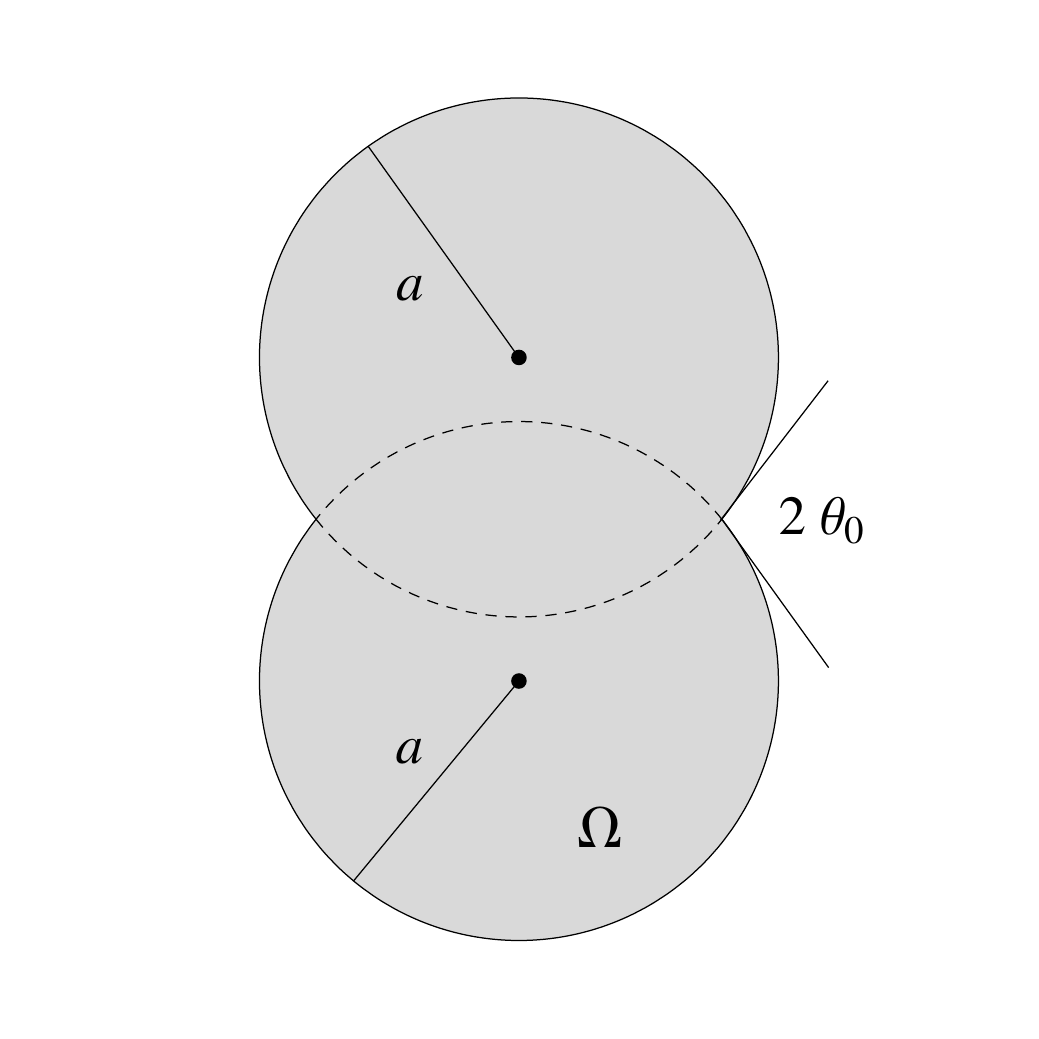}
      \caption{}
    \end{subfigure}\\
    \begin{subfigure}{0.48\textwidth}
      \centering
      \includegraphics[width=6.3cm,height=4.8cm]{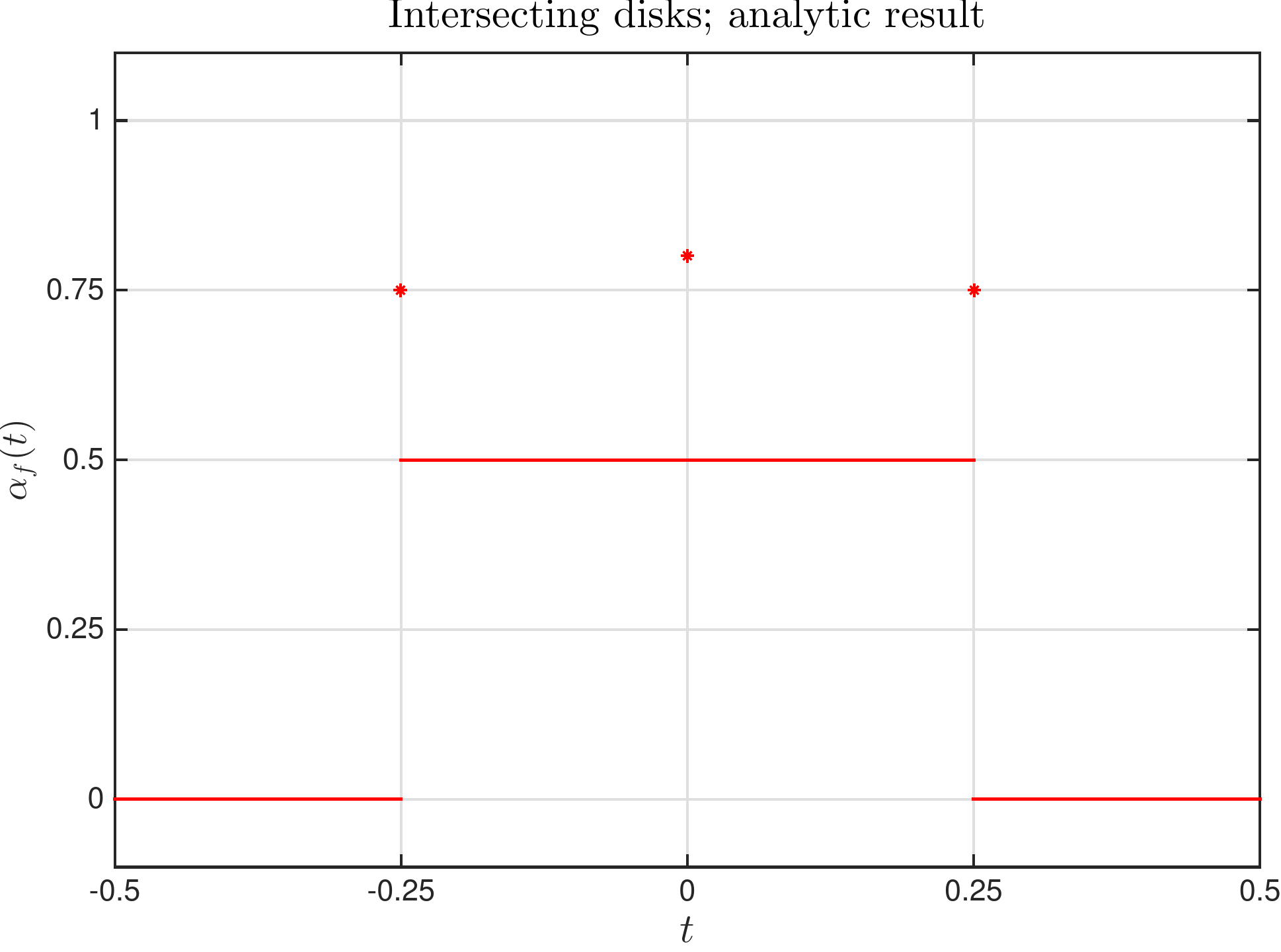}
      \caption{}
    \end{subfigure}%
    \begin{subfigure}{0.48\textwidth}
      \centering
      \includegraphics[width=6.3cm,height=4.8cm]{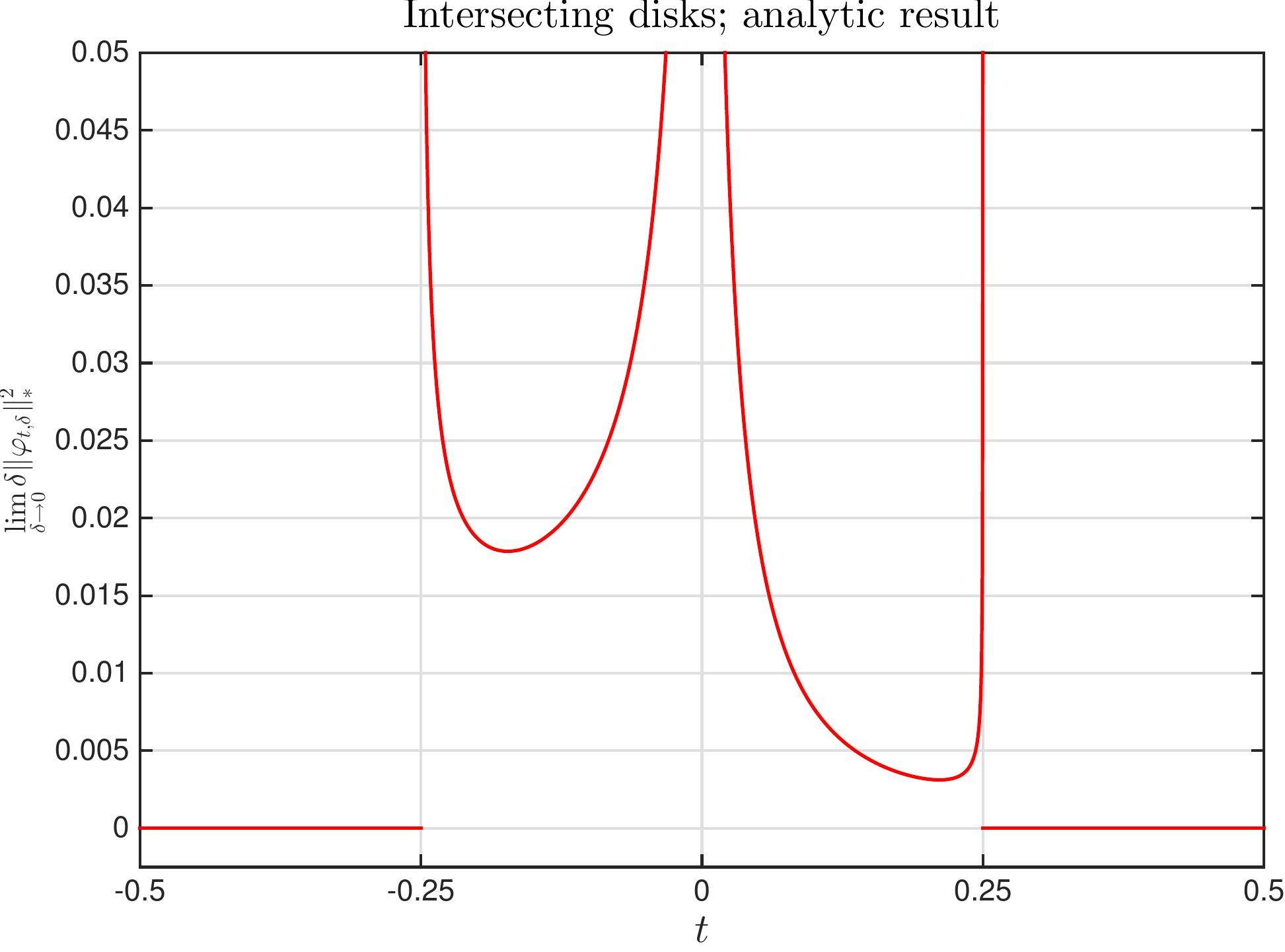}
      \caption{}
    \end{subfigure}
    \begin{subfigure}{0.48\textwidth}
      \centering
      \includegraphics[width=6.5cm,height=4.8cm]{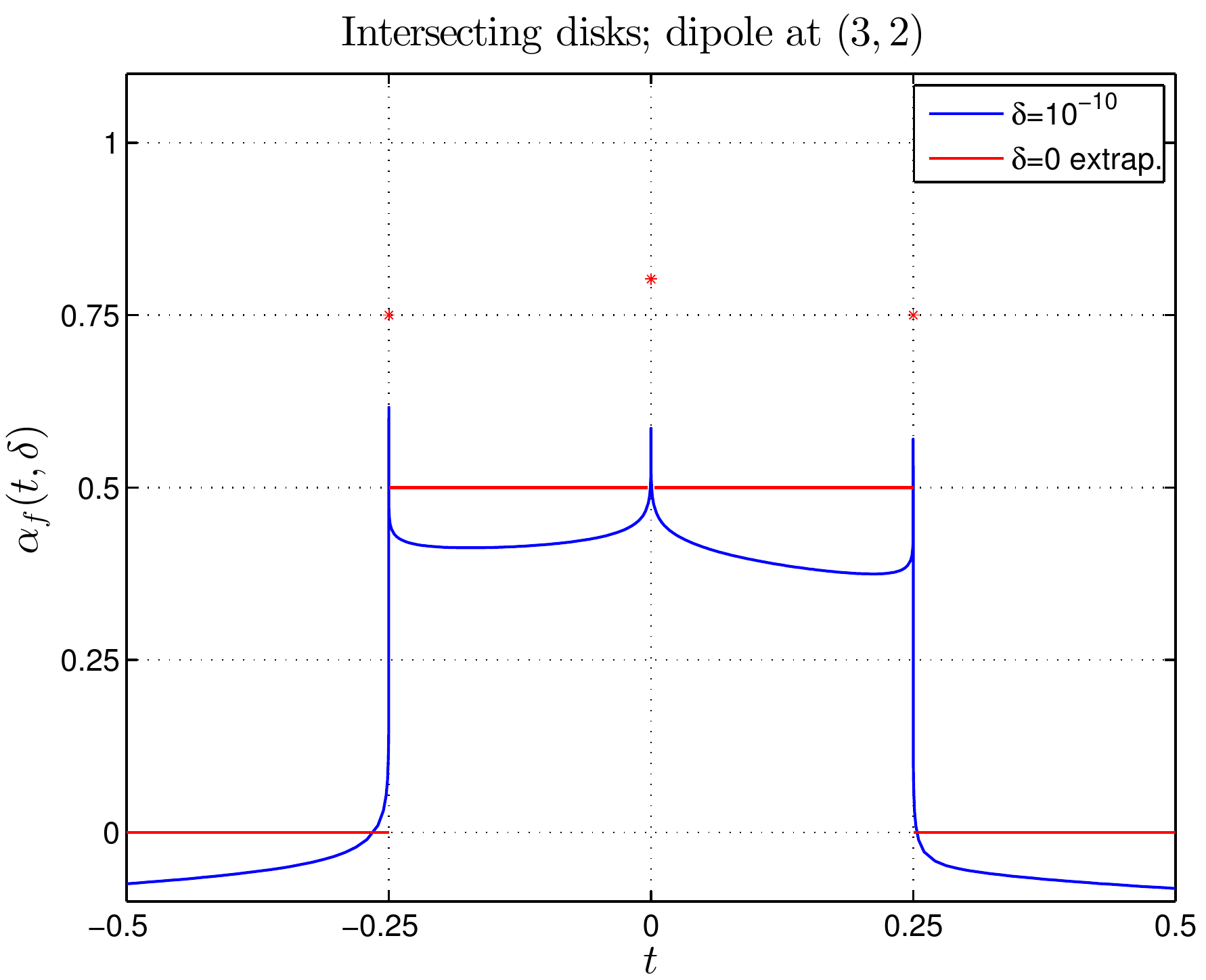}
      \caption{}
    \end{subfigure}%
    \begin{subfigure}{0.48\textwidth}
      \centering
      \includegraphics[width=6.5cm,height=4.8cm]{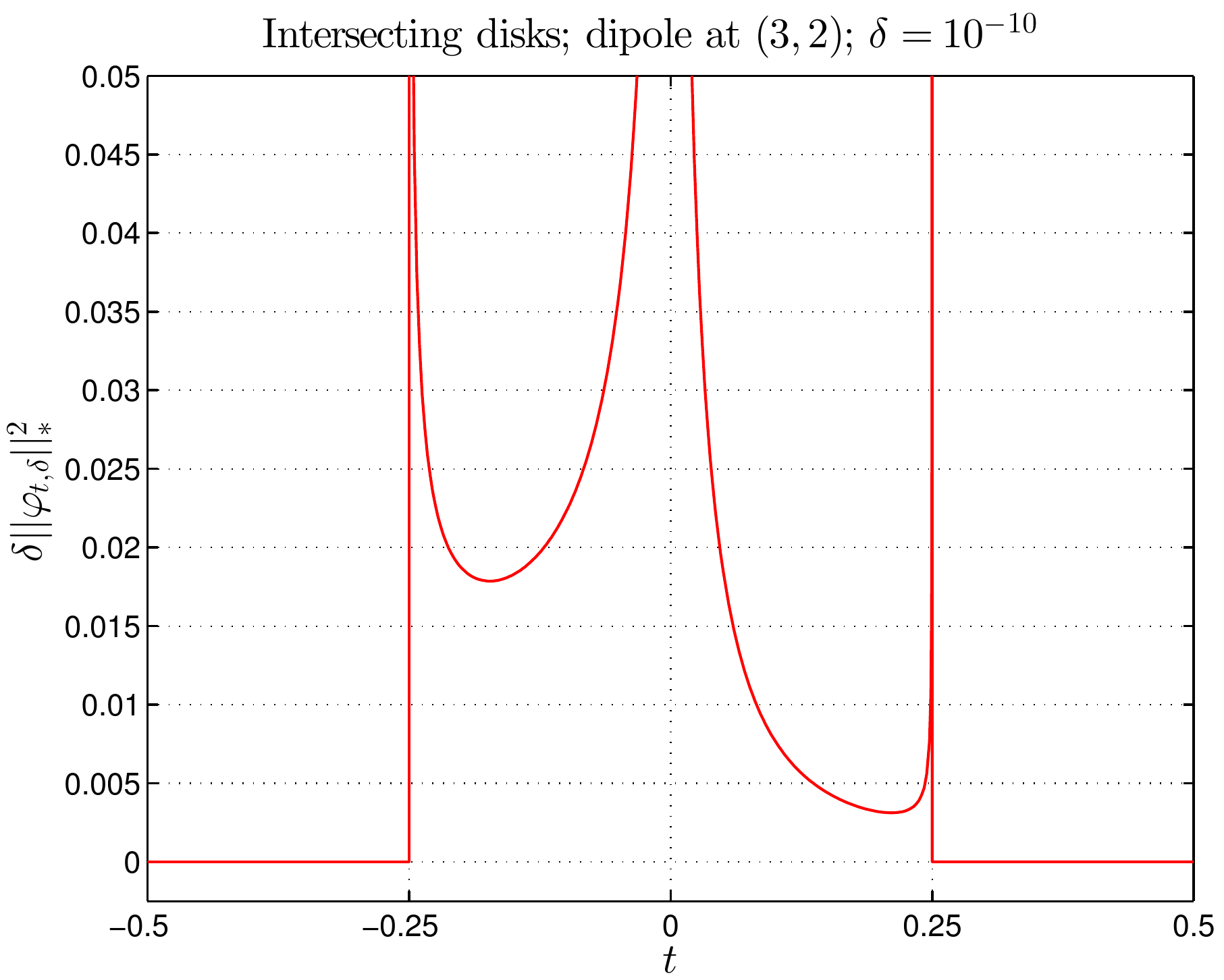}
      \caption{}
    \end{subfigure}
      \begin{subfigure}{0.48\textwidth}
      \centering
      \includegraphics[width=6.5cm,height=4.8cm]{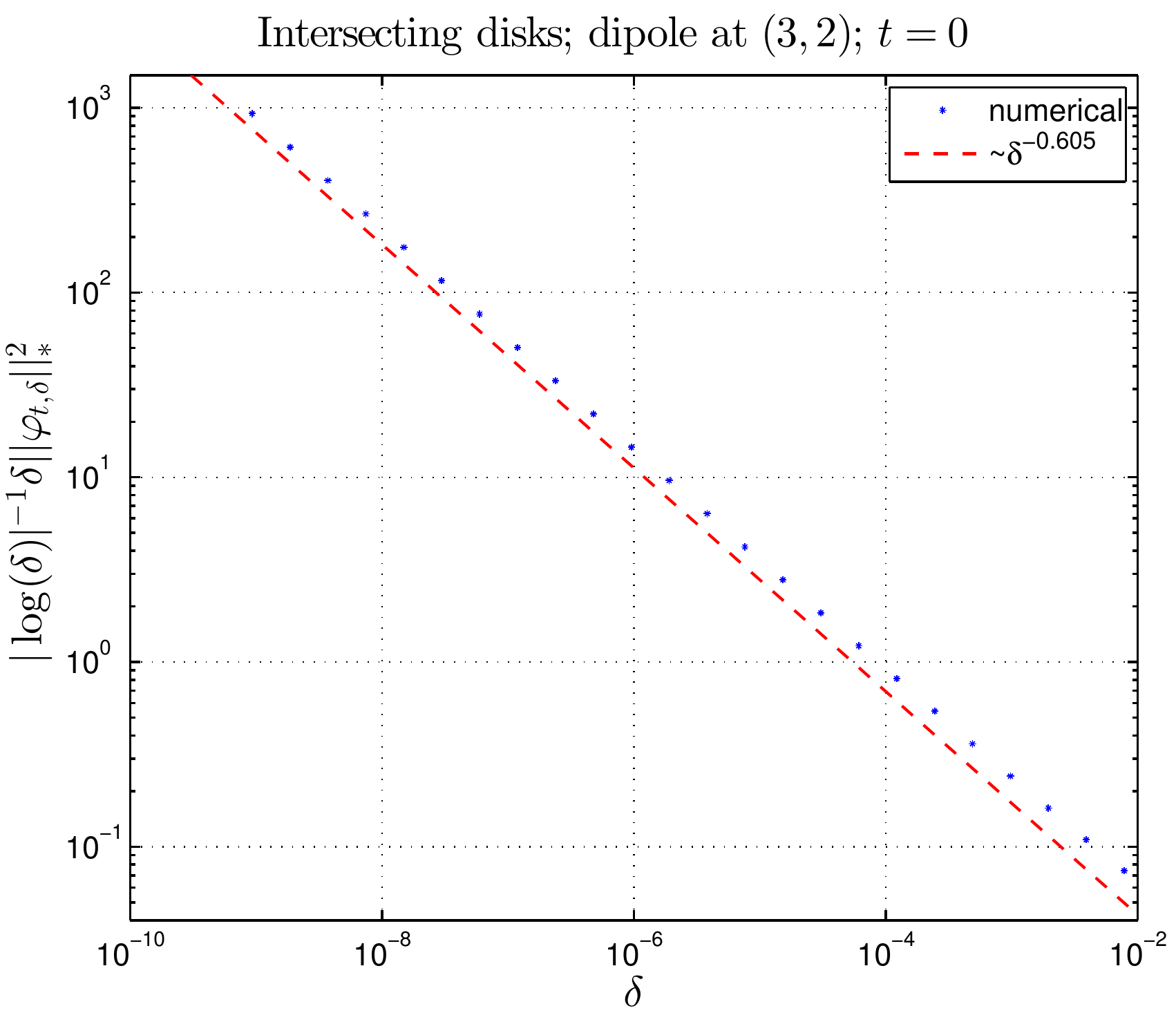}
      \caption{}
    \end{subfigure}%
    \begin{subfigure}{0.48\textwidth}
      \centering
      \includegraphics[width=6.5cm,height=4.8cm]{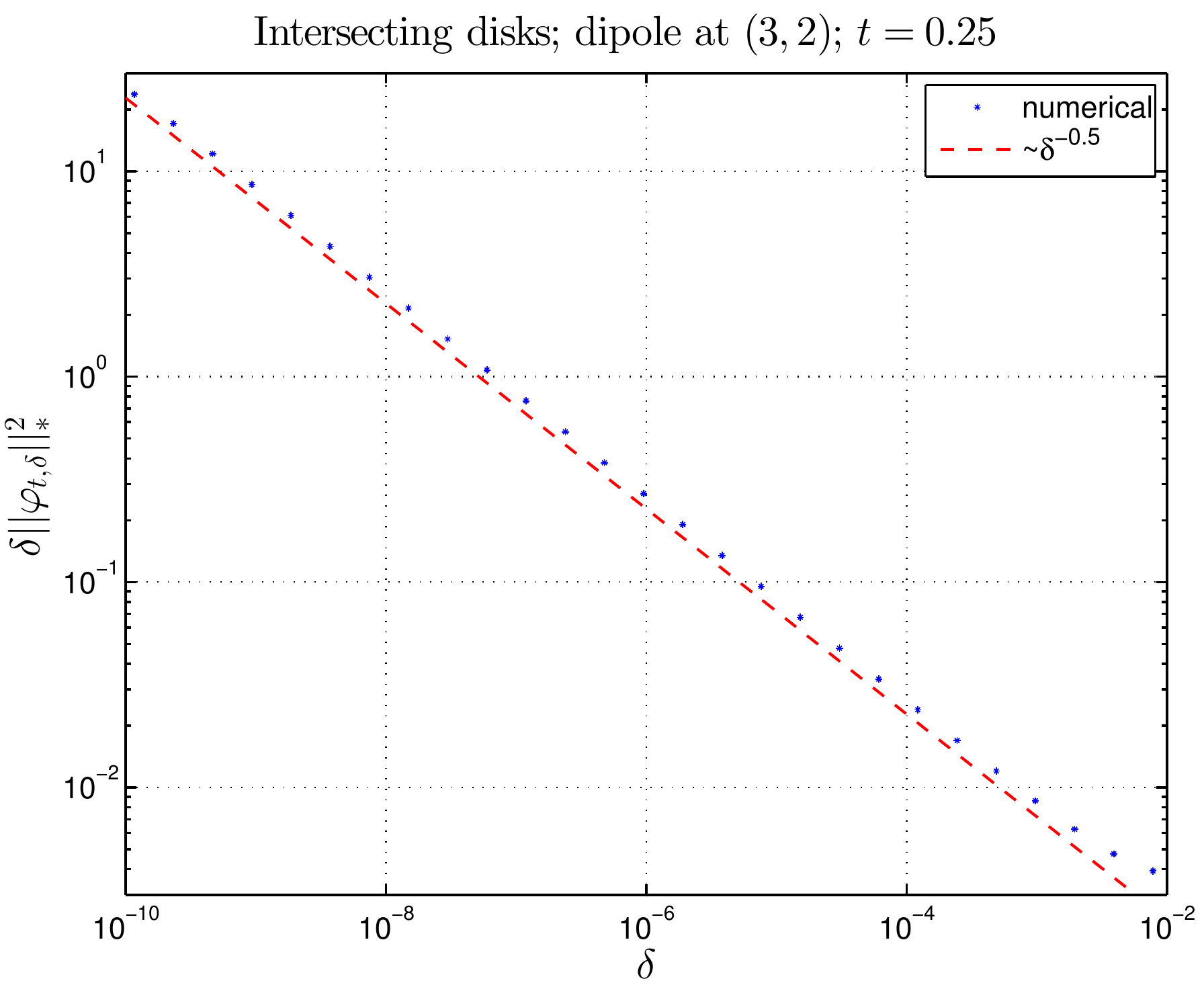}
      \caption{}
    \end{subfigure}

\caption{Spectrum of the intersecting disks taken by one dipole source. (a) illustrates the intersecting disks (in blue) and the source (location as a dot and orientation as an arrow). (b) illustrates the intersecting disks whose exterior angle, say $2\theta_0$, at the two corners is smaller than $\pi$.
(c) and (d) are analytical values; (e), (f), (g), and (h) are numerical results
using the RCIP-accelerated Nystr{\"o}m solver.}
\label{fig:ID;onedipole}
\end{figure}

Numerical results for the spectrum of the NP operator, obtained with
the RCIP-accelerated Nystr{\"o}m solver, are now compared with
analytical results. The domain $\Omega$ is that of the intersecting
disks in Fig.\;\ref{fig:ID;onedipole}(b). We let the disk radius be
$a=2$ and the exterior angle at the two corners be
$\theta_0={\pi}/{4}$.

Let us briefly review the analytical results for intersecting disks
obtained in \cite{KLY}. The spectrum of the NP operator has the range
$[-b,b]$ with $b=0.25$, which is in agreement with \eqnref{bounds},
and consists only of absolutely continuous spectrum. Let
$(\Psi_1(z),\Psi_2(z))$ be the bipolar coordinates with two foci $(\pm
c,0)$ located at the corners of the intersecting disks, {\it i.e.},
\begin{displaymath}
\Psi_1(z)+i\Psi_2(z)=\mbox{Log}
\left(\frac{z+c}{z-c}\right),\quad
c=a\sin\theta_0\,,
\end{displaymath}
where $\mbox{Log}$ is the logarithm with the principal branch. For the
dipole field $f_z$, oriented in a suitable direction, we have
\begin{equation}
\label{eqn:alphafz}
\alpha_{f_z}(t)=\begin{cases}
0\,,\quad&|t|>b\,,\\
0.5\,,\quad& 0<|t|<b\,,\\
0.5\left(1+\frac{|\Psi_2(z)|}{\theta_0}\right),\quad&t=0\,,\\
0.75\,,\quad&t=\pm b\,.
\end{cases}
\end{equation}
More precisely, $\delta\|\varphi_{t,\delta}\|_*^2$ converges to a
positive number as $\delta\to 0$ for $0<|t|<b$ with a limit that can
be expressed in terms of elementary functions of bipolar coordinates.
Furthermore, $\delta^{3/2}\|\varphi_{t,\delta}\|_*^2$ converges at
$|t|=b$ and
\begin{equation}
\label{eqn:phizero}
\left|\log
  \delta\right|^{-1}\delta^{1+{|\Psi_2(z)|}/{\theta_0}}\|\varphi_{t,\delta}\|_*^2\quad\mbox{converges
  at }t=0\,.
\end{equation}
In view of \eqnref{eqn:alphafz}, the indicator function
$\alpha_{f_z}(0)$ depends on the location of the dipole source. It
increases as $z$ approaches $\p\Omega$, but never reaches one.

\smallskip

\noindent{\textbf{One dipole field}.} We first consider the spectrum of
the intersecting disks as excited by a single dipole source located at
$z=(3,2)$ and orientated in the direction $d=(1,1)/\sqrt{2}$. In
Fig.\;\ref{fig:ID;onedipole}, images (c) and (d) show analytical
values of $\alpha_{f_z}(t)$ from \cite{KLY} while (e), (f), (g), and (h) show numerical
results. Table\;\ref{table:analyticRCIP} compares values from
Fig.\;\ref{fig:ID;onedipole}(d,f) and shows that the numerical results
of (f) for $\delta=10^{-10}$ exhibit a $6$-digit agreement with the
analytical results of (d) for
$\lim_{\delta\rightarrow0}\delta\|\varphi_{t,\delta}\|_*^2$ when $t$
is away from $0,\pm b$. In Fig.\;\ref{fig:ID;onedipole}(e), numerical
values for $\alpha_f(t,\delta)$ (in blue) are computed with
$\delta=10^{-10}$, and those for $\alpha_f(t)$ (in red) are
extrapolated from the limit behavior of
$\delta\|\varphi_{t,\delta}\|_*^2$ for $|t|\neq 0$. At $t=0$,
$\left|\log
  \delta\right|^{-1}\delta^{1+{|\Psi_2(z)|}/{\theta_0}}\|\varphi_{t,\delta}\|_*^2$
is instead considered because of the property \eqnref{eqn:phizero}.
Note that the $\log\delta$ factor does not affect the limit value of
$\alpha_f(t,\delta)$ as $\delta\to 0$. We conclude that analytical
values and numerical results of $\alpha_f(t)$, red graphs in
Fig.\;\ref{fig:ID;onedipole}(c,e), coincide for $t\neq0$ and have very
similar values at $t=0$.

\begin{table}[ht]
\centering 
\begin{tabular}{|c|| l | l |} 
\hline 
$t$&$\lim_{\delta\rightarrow0}\delta\|\varphi_{t,\delta}\|_*^2$ (analytical)&$\delta\|\varphi_{t,\delta}\|_*^2$, $\delta=10^{-10}$ (numerical)\\[0.5ex]
\hline 
$-$0.3&0&2$\cdot 10^{-11}$\\[0.3ex]
$-$0.2& 0.018710399304385&0.0187104\\[0.3ex]
 $-$0.1&  0.022245420816273&0.0222454\\ [0.3ex]
 $+$0.1&0.007687535353992&0.00768753\\[0.3ex]
  $+$0.2&0.003180101918936&0.00318010\\[0.3ex]
  $+$0.3&0&8$\cdot 10^{-12}$\\
\hline 
\end{tabular}
\caption{Spectrum of the intersecting disks taken by one dipole source located at $(3,2)$.}
\label{table:analyticRCIP}
\end{table}

\smallskip

\begin{figure}[htbp!]
\centering\includegraphics[width=5.8cm]{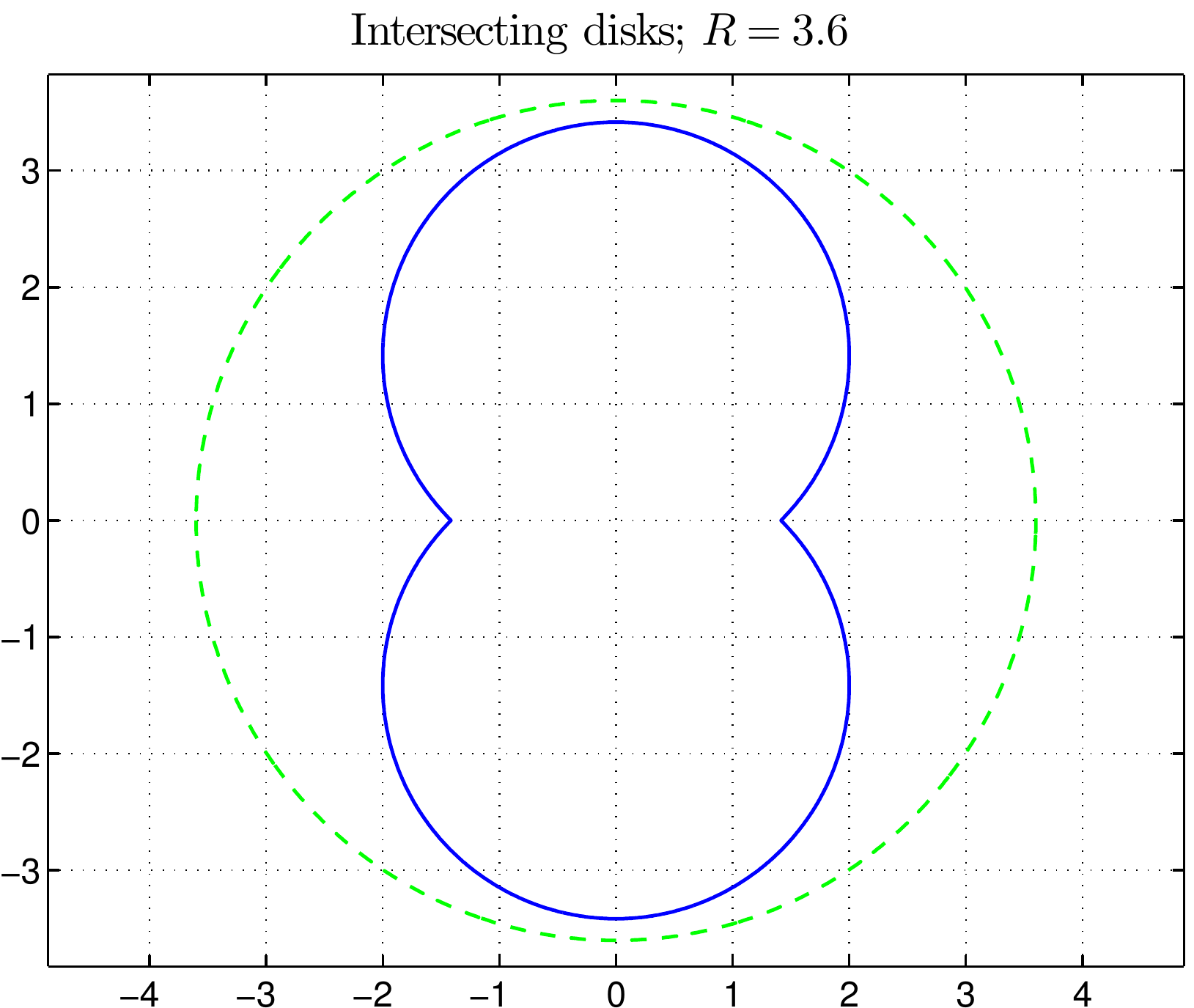}\hskip 7mm
\centering\includegraphics[width=6cm, height=5cm]{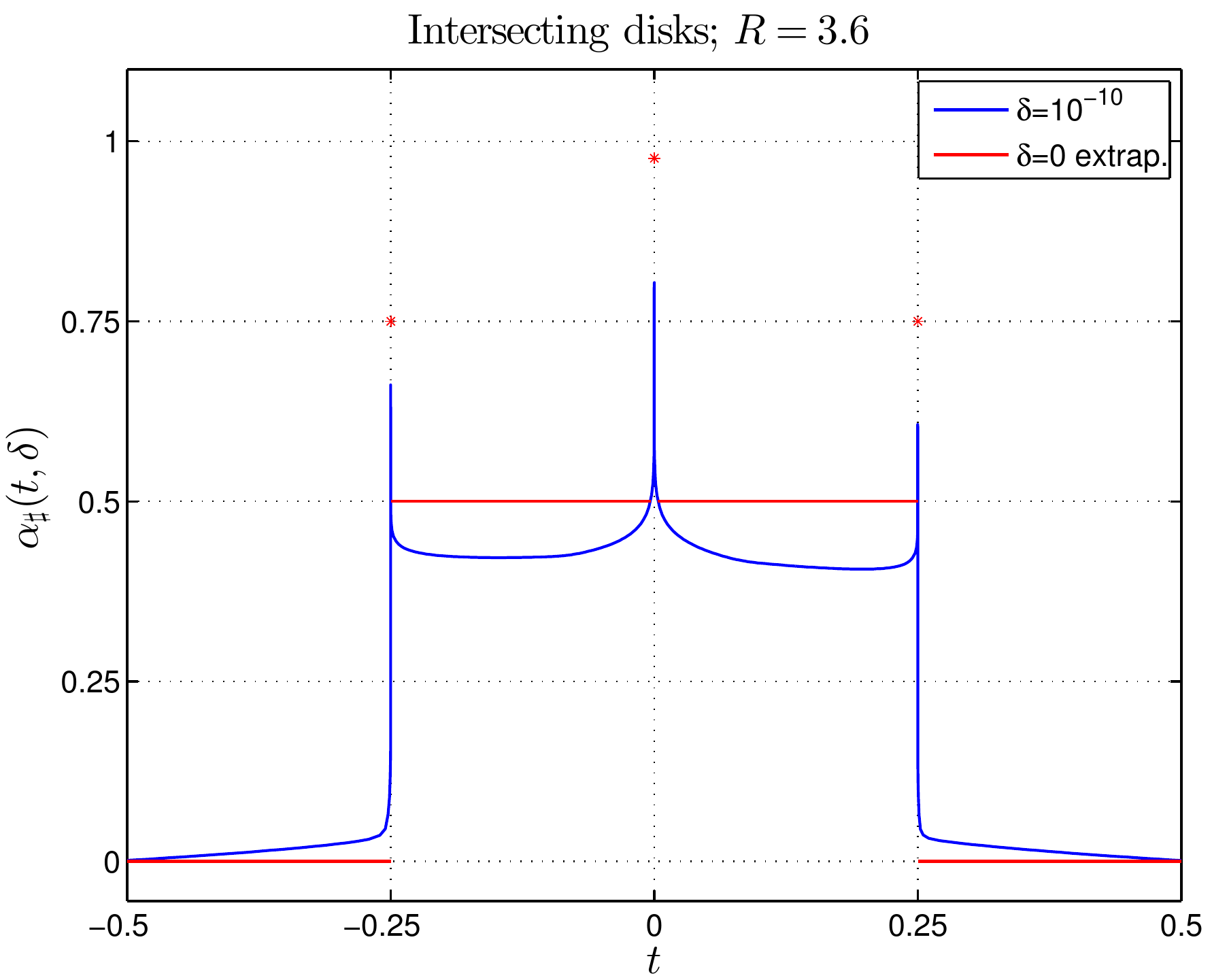}\\[2mm]
\centering\includegraphics[width=6cm, height=5.cm]{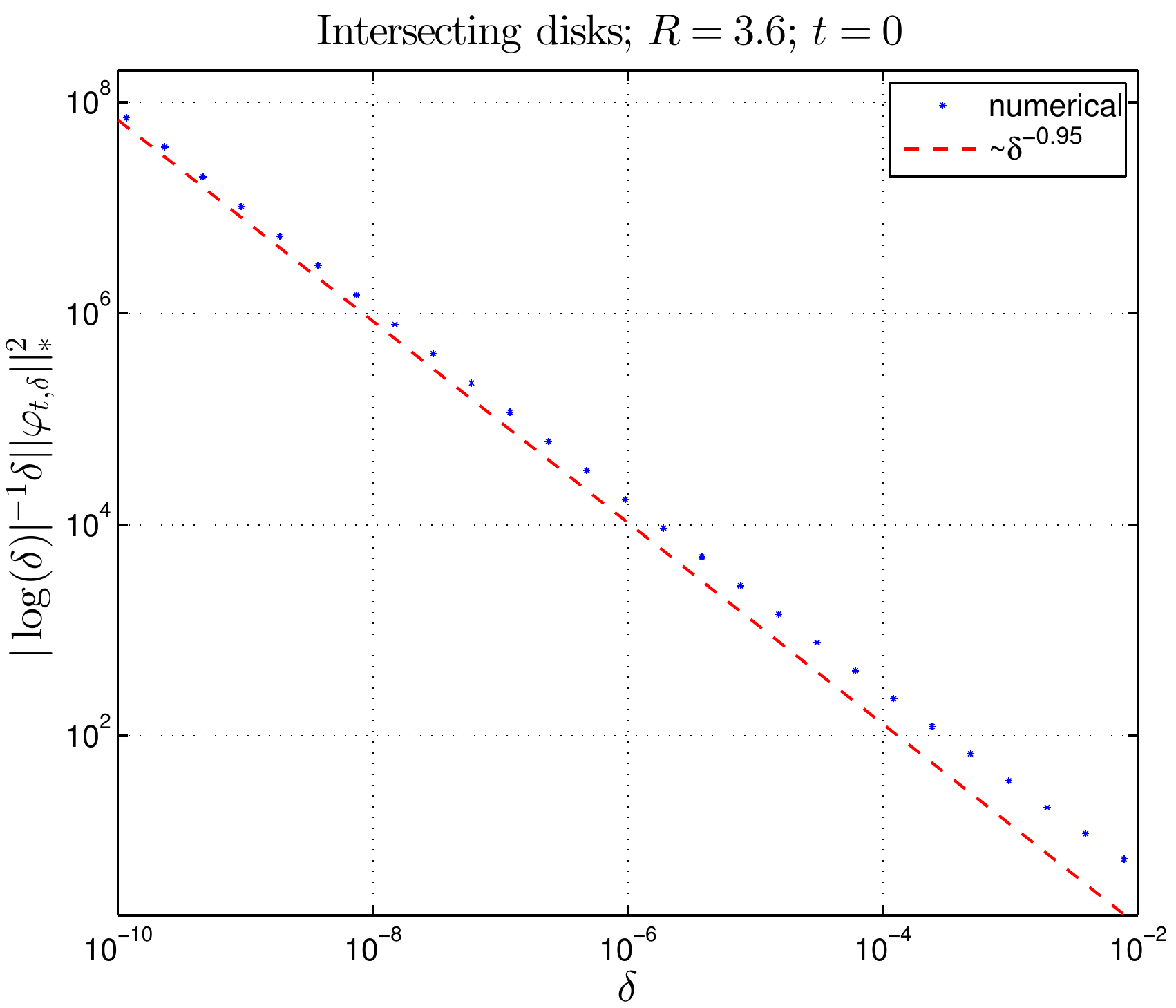}\hskip 7mm
\centering\includegraphics[width=6cm, height=5.cm]{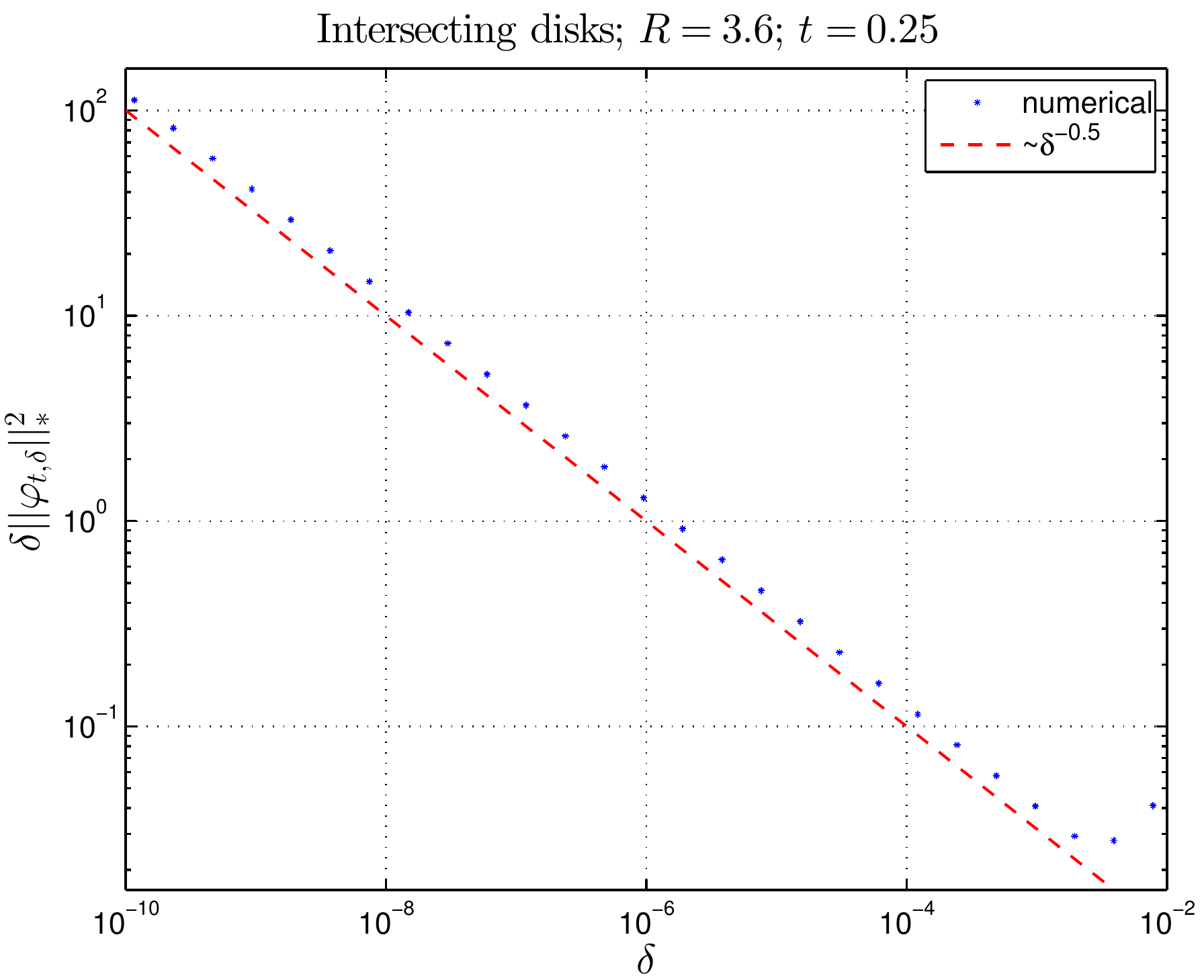}
\caption{Spectrum of the intersecting disks taken over many dipole fields. Dipole sources are situated on an enclosing circle.}
\label{fig:ID2}
\end{figure}

\noindent{\textbf{A large number of dipole fields}.}
Fig.\;\ref{fig:ID2} shows $\alpha_\sharp(t)$ taken over dipole sources
situated on the enclosing circle with radius $R=3.6$. The obtained
extrapolated $\alpha_\sharp(t)$ coincides with $\alpha_f(t)$ in
Fig.\;\ref{fig:ID;onedipole}(e) for $t\neq0$. At $t=0$, the indicator
function $\alpha_\sharp(t)$ achieves a larger value than $\alpha_f(t)$
since it involves the maximum over multiple sources and some of sources attain bigger values in $|\Psi_2(z)|$.

\subsection{Triangle}\label{subsec:tri}

\begin{figure}[htbp]
\centering\includegraphics[width=6cm]{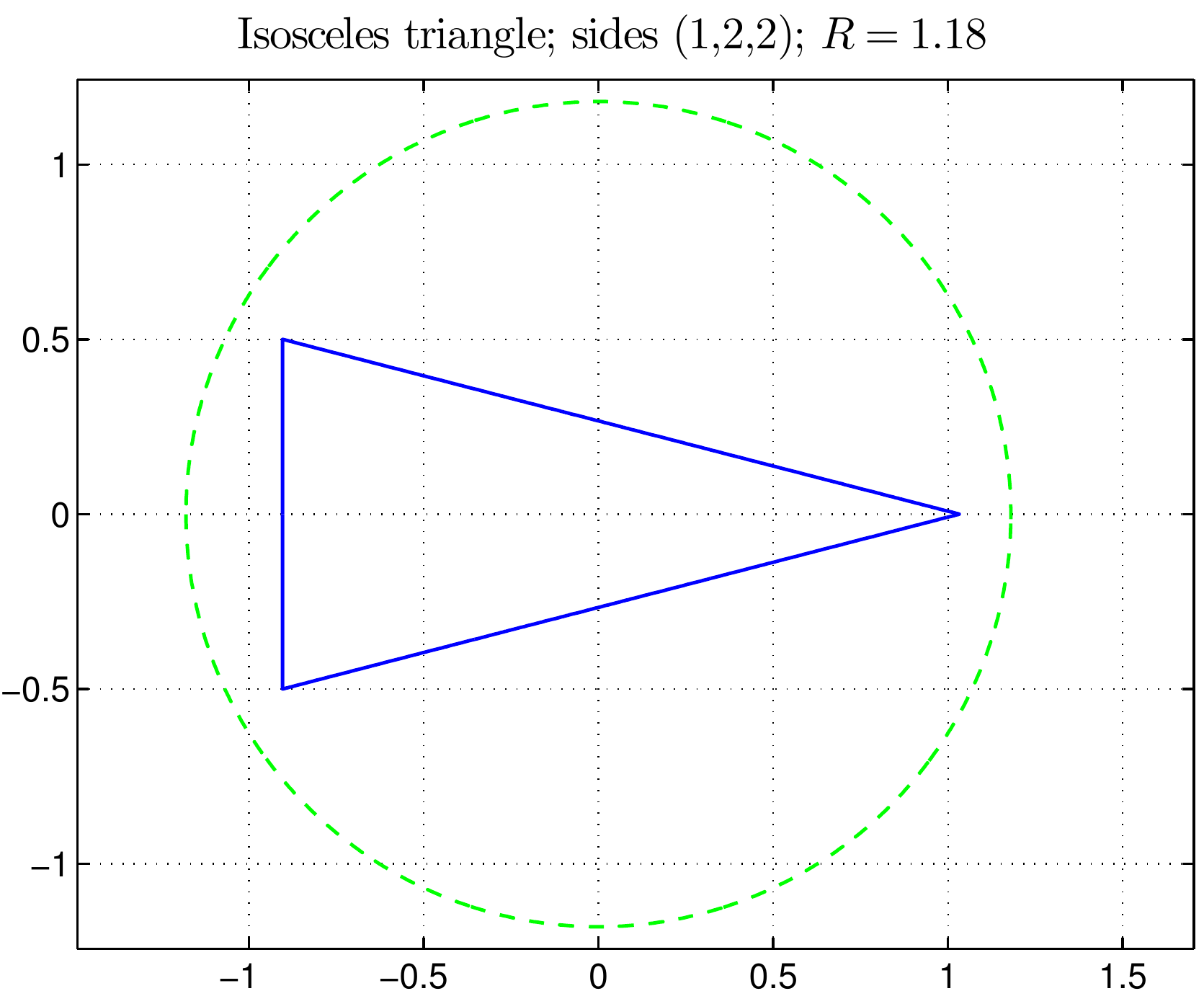}\hskip 7mm
\centering\includegraphics[width=6cm]{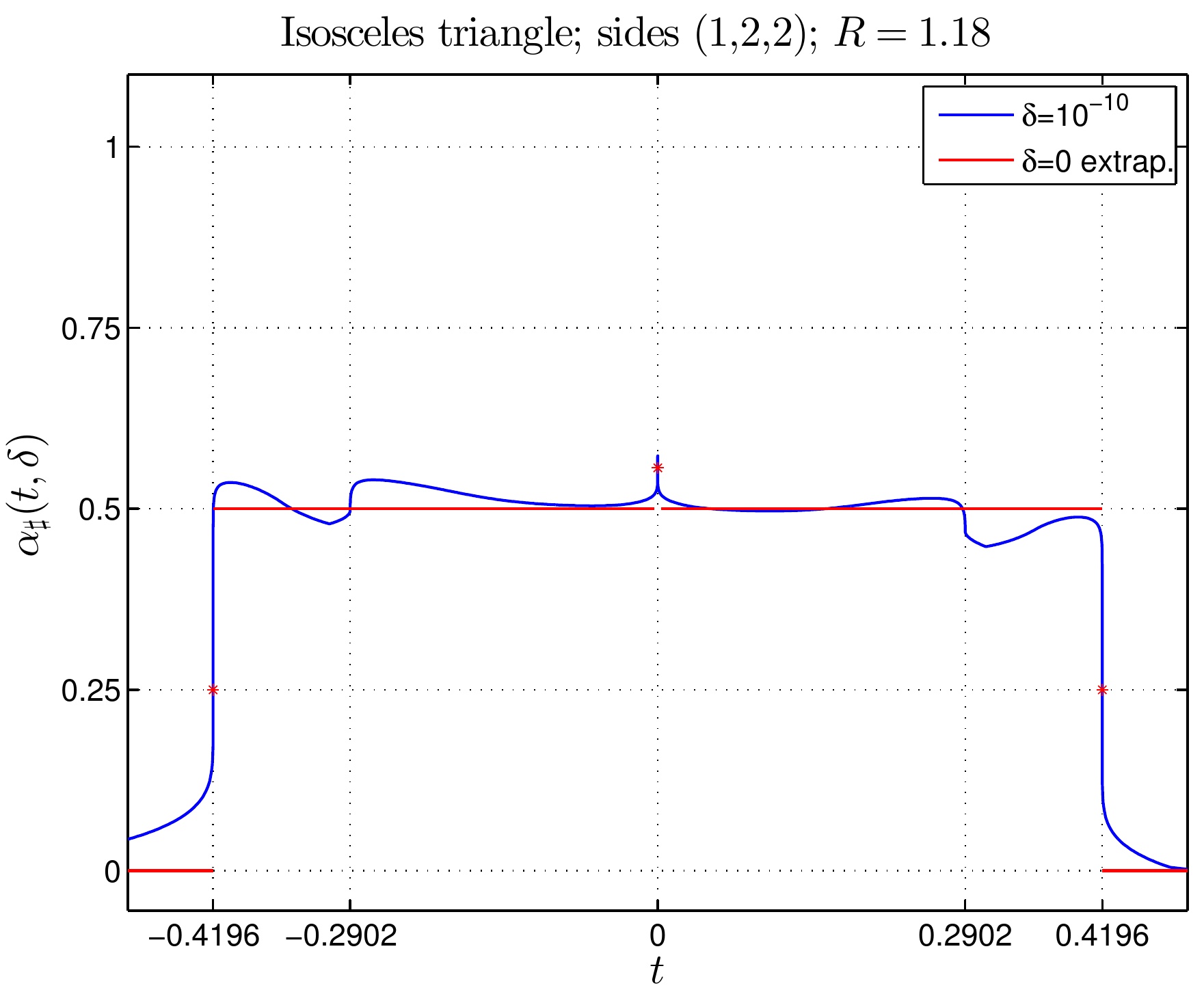}\\[2mm]
\centering\includegraphics[width=6cm]{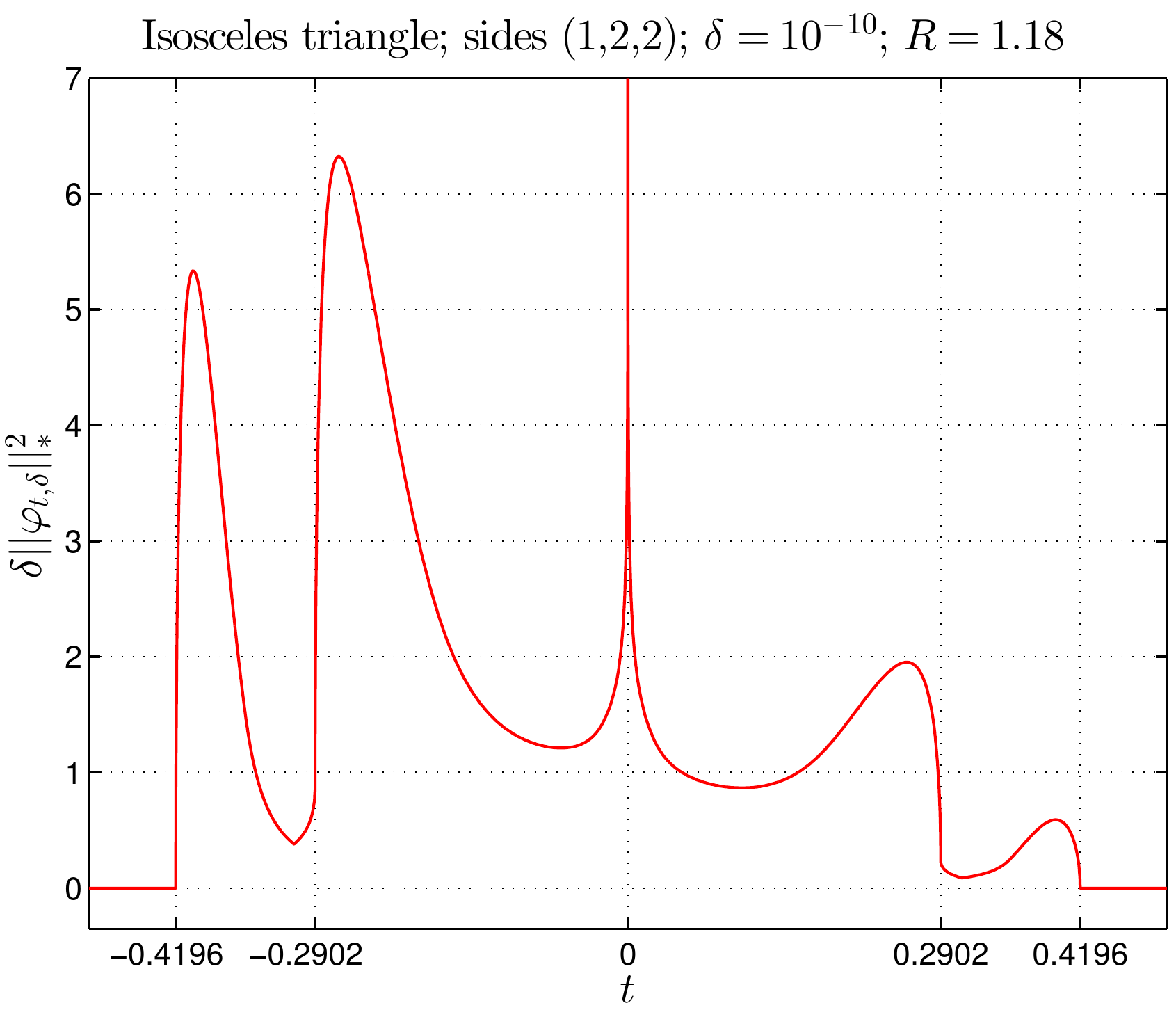}\hskip7mm
\centering\includegraphics[width=6cm]{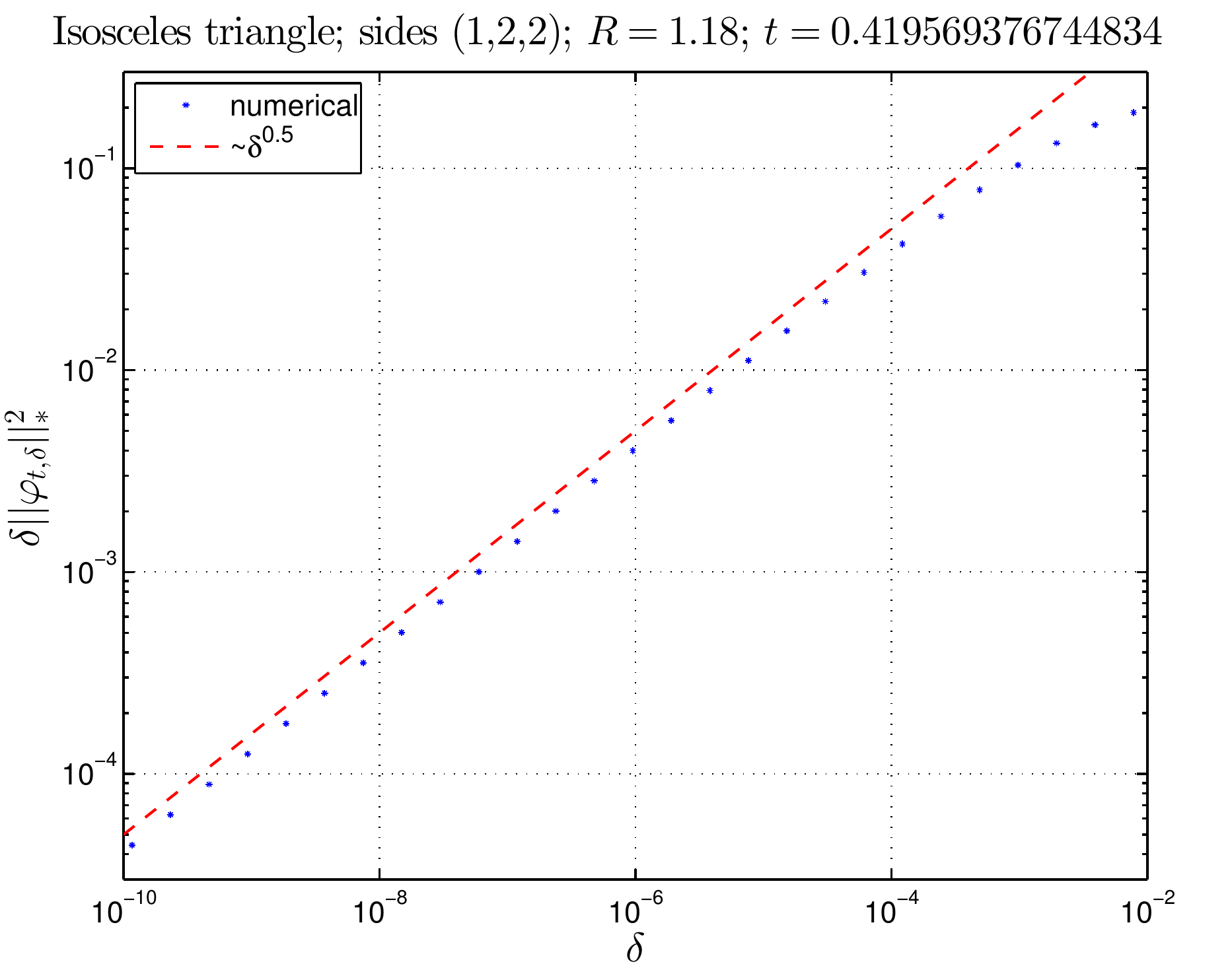}
\caption{Spectrum of the isosceles triangle with sides $1$, $2$ and $2$. The values of $0.5(1-\theta/\pi)$ for interior angles, say $\theta$, are approximately $0.4196$ and $0.2902$. The larger number $0.4196$ bounds the essential spectrum. While the indicator function $\alpha_\sharp(t)$ changes only at zero and $0.4196$, the functions $\alpha_\sharp(t,\delta)$ and $\delta\|\varphi_{t,\delta}\|_*^2$ for $\delta=10^{-10}$ show dynamic changes near $0.2902$ as well.}
\label{fig:IT}
\end{figure}

Fig.\;\ref{fig:IT} show the spectrum of an isosceles triangle. The interval of continuous spectrum is determined by the smallest interior angle according to \eqnref{bounds}.

\subsection{Rectangles and superellipses}\label{subsec:rec}

\begin{figure}[htbp!]
\begin{center}
      \includegraphics[width=6.4cm]{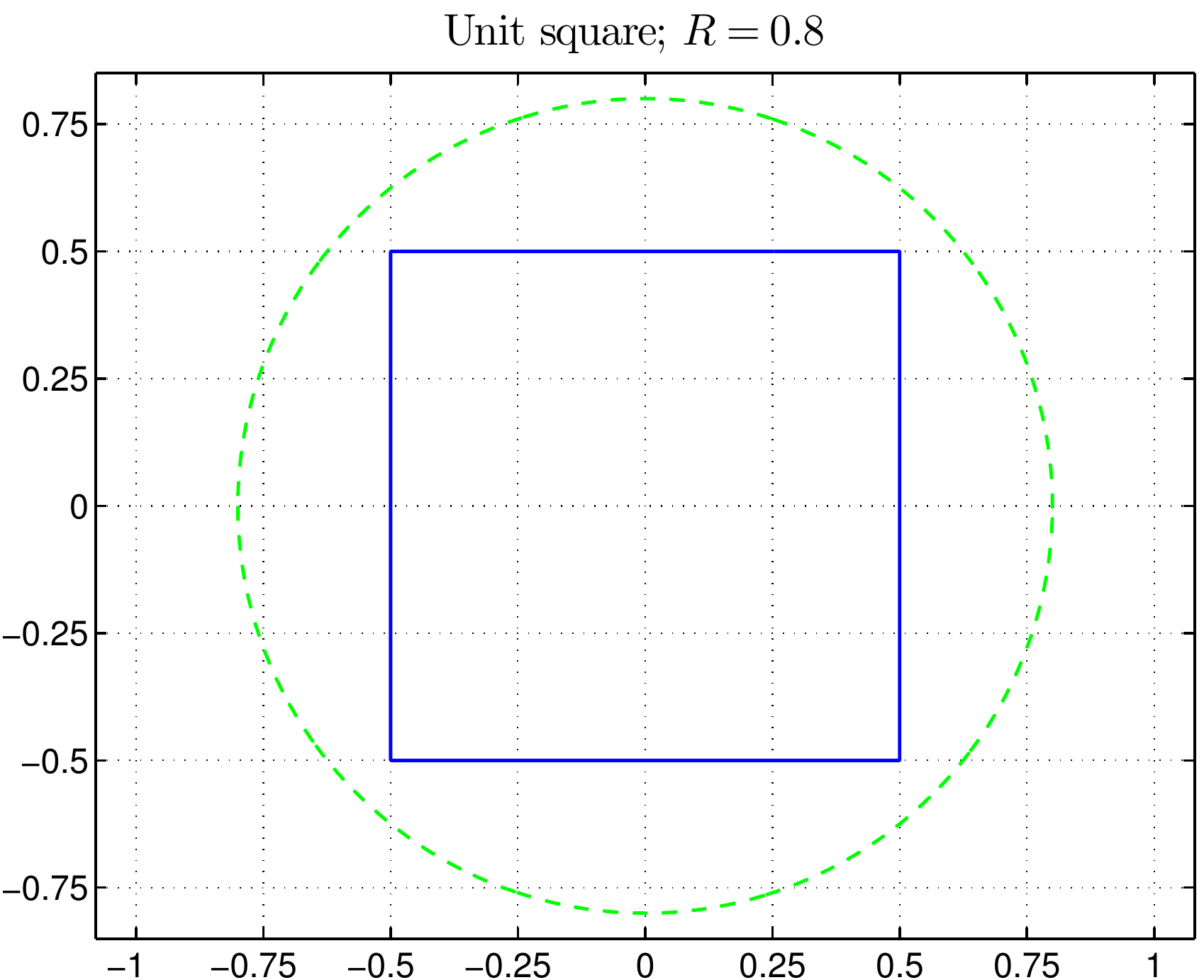}\hskip 7mm
      \includegraphics[width=6.5cm]{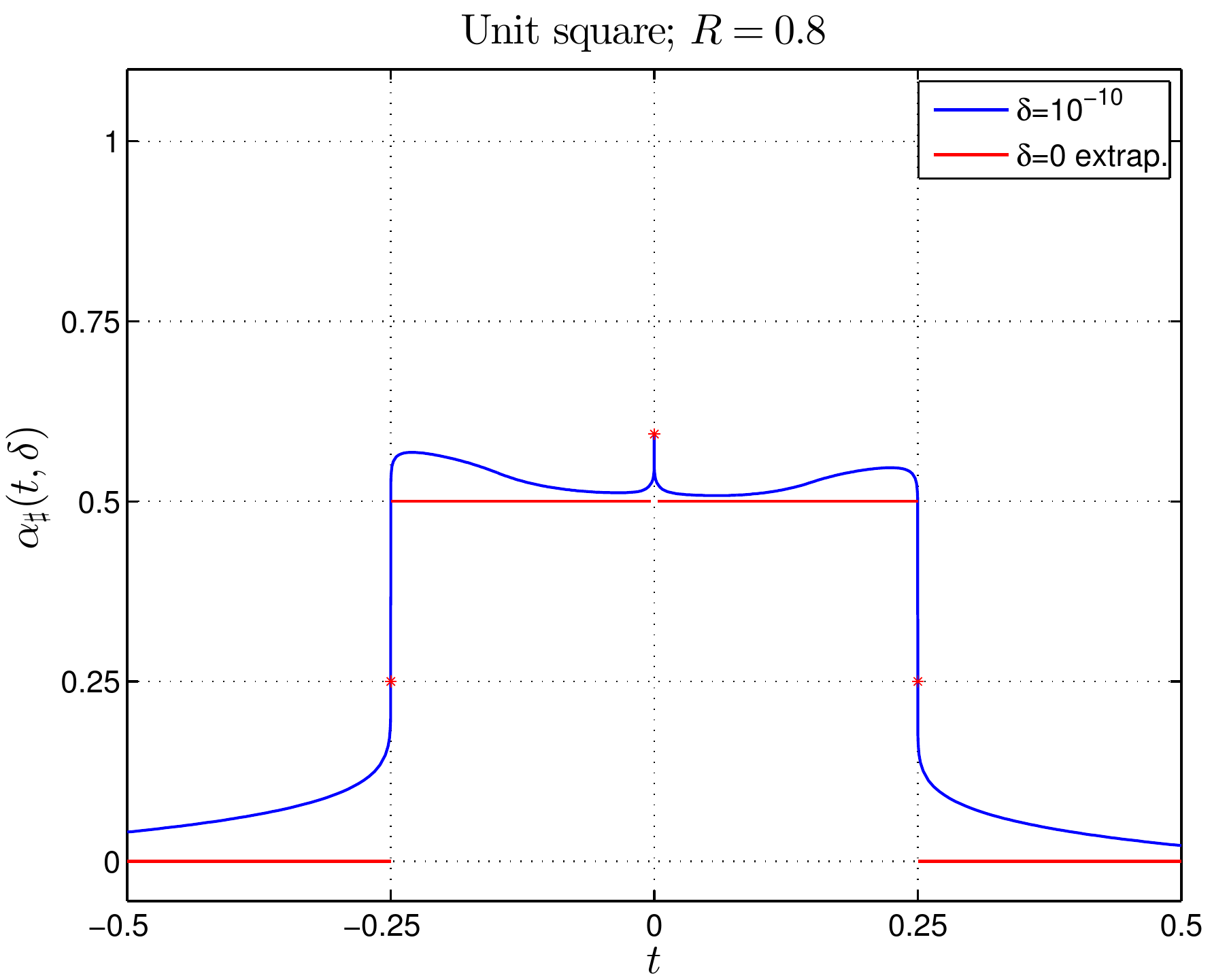}\\[2mm]
       \includegraphics[width=6.4cm]{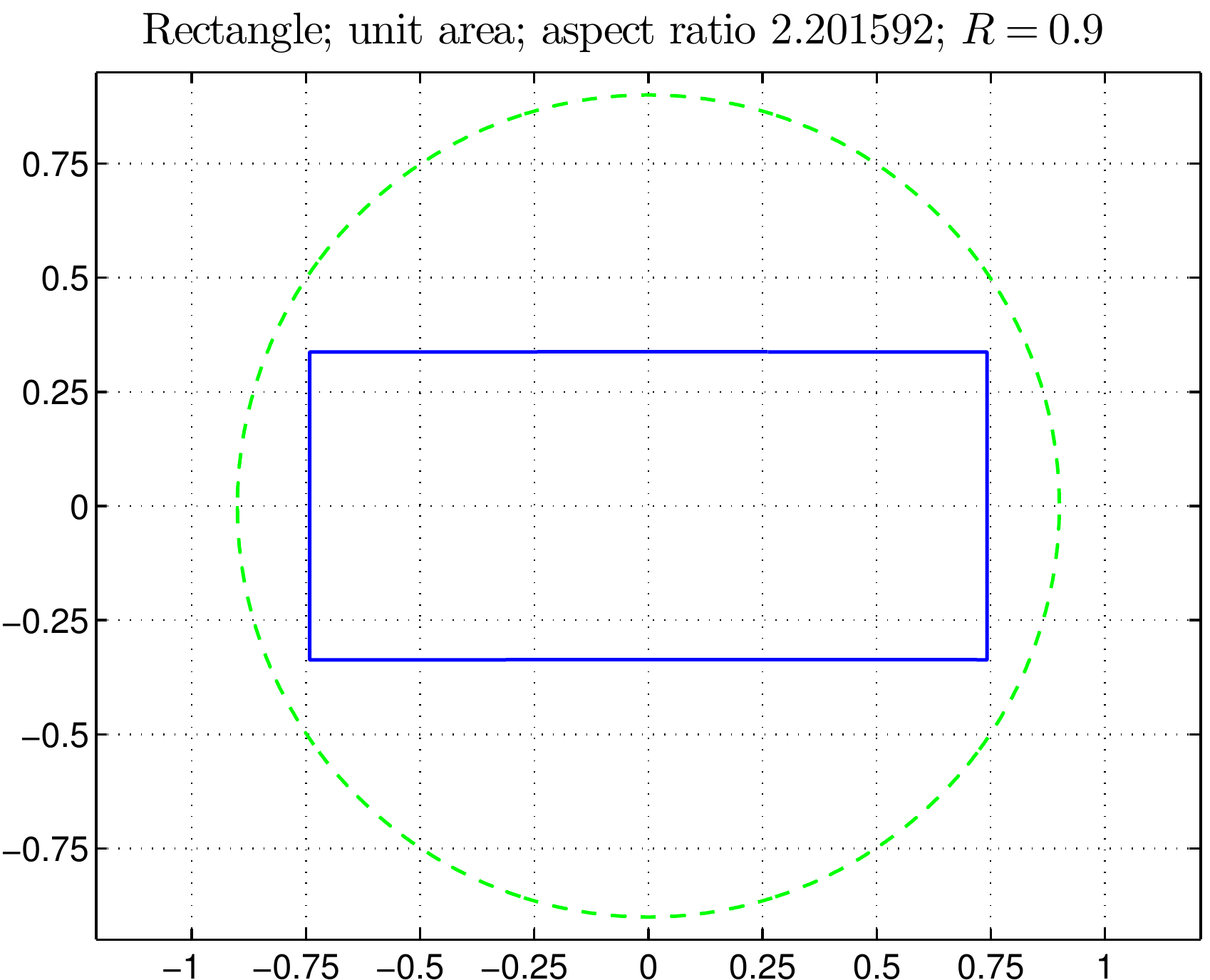}\hskip 7mm
         \includegraphics[width=6.5cm]{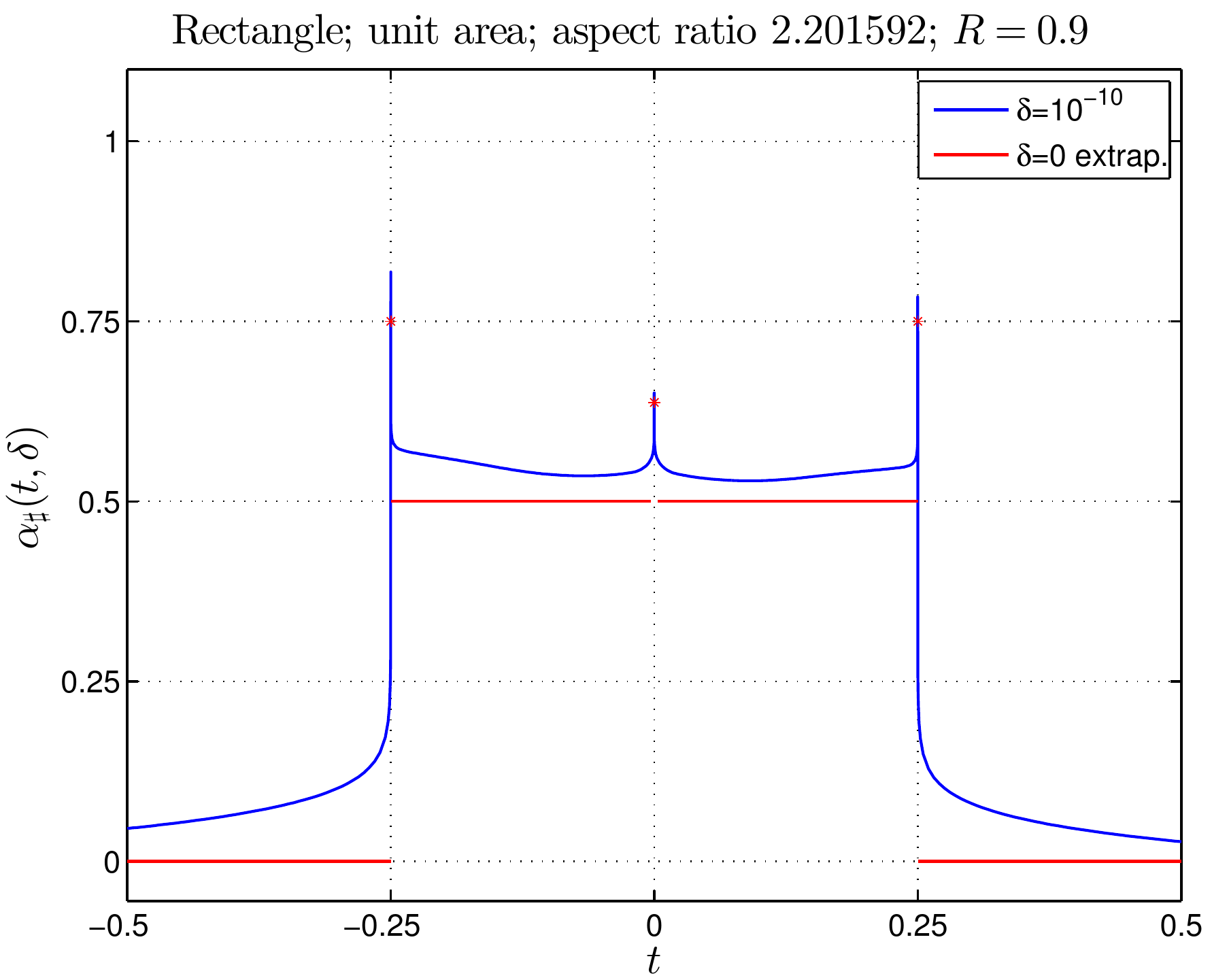}  \\[2mm]
               \includegraphics[width=6.4cm]{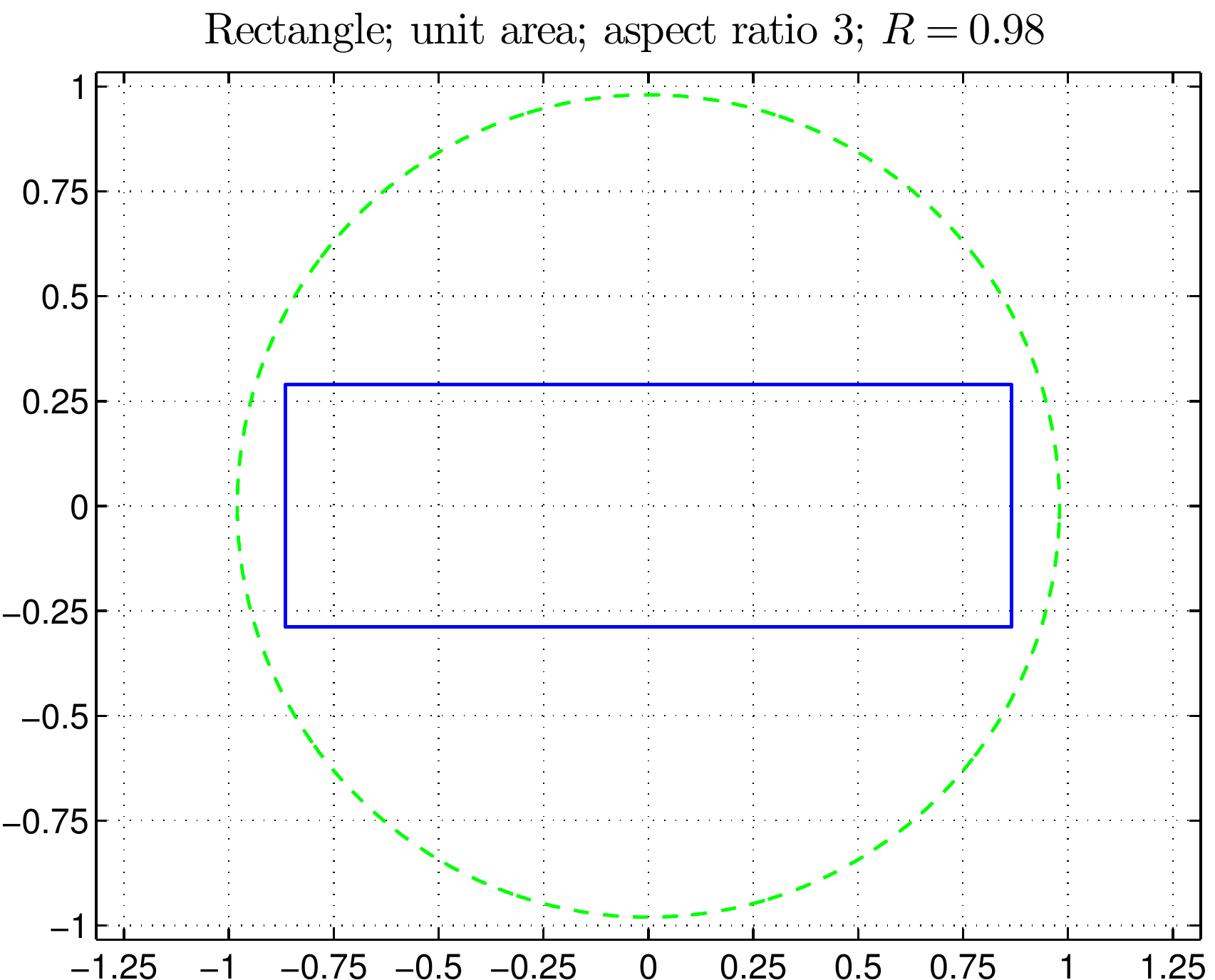}\hskip 7mm
        \includegraphics[width=6.5cm]{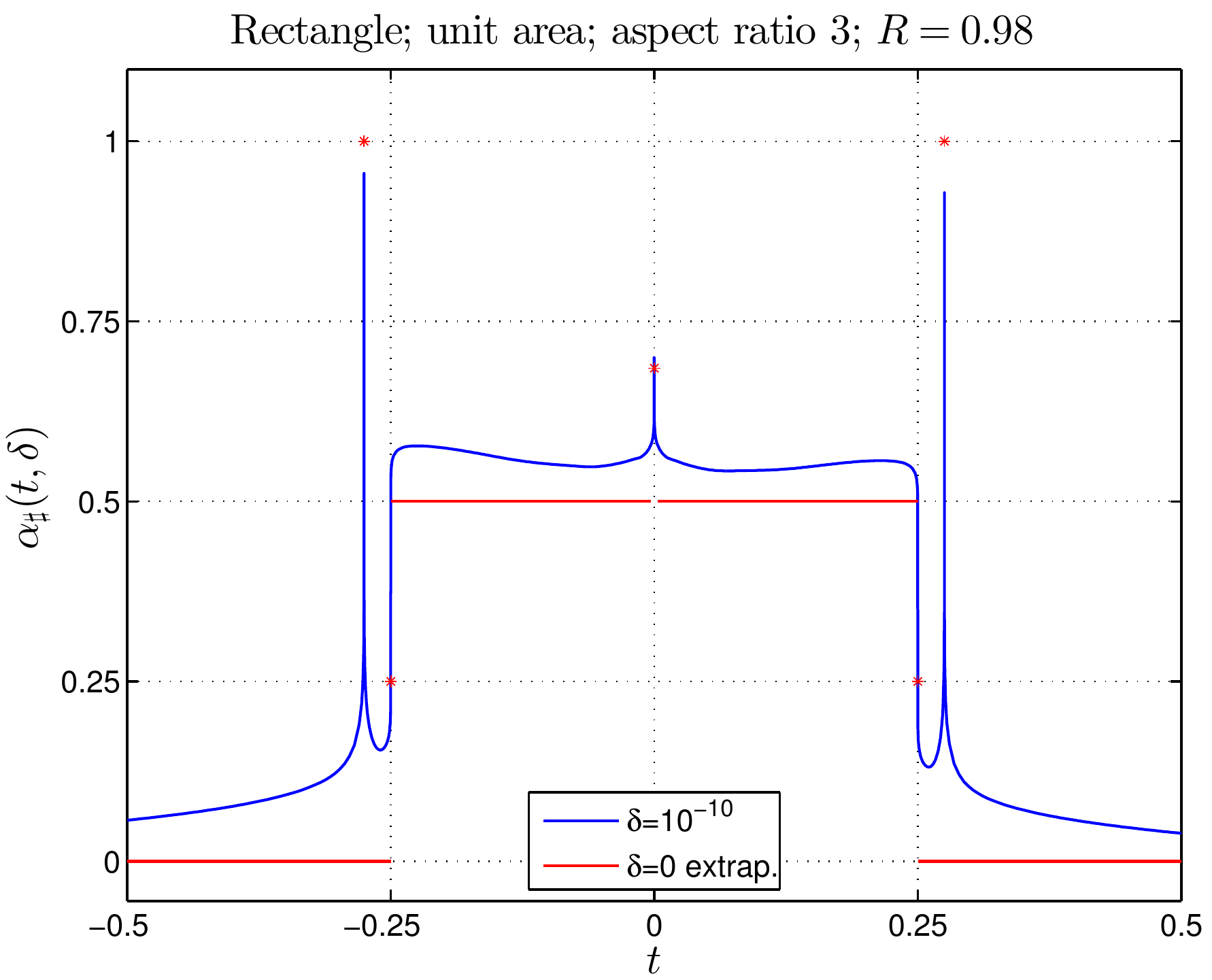} \\[2mm]
          \includegraphics[width=6.4cm]{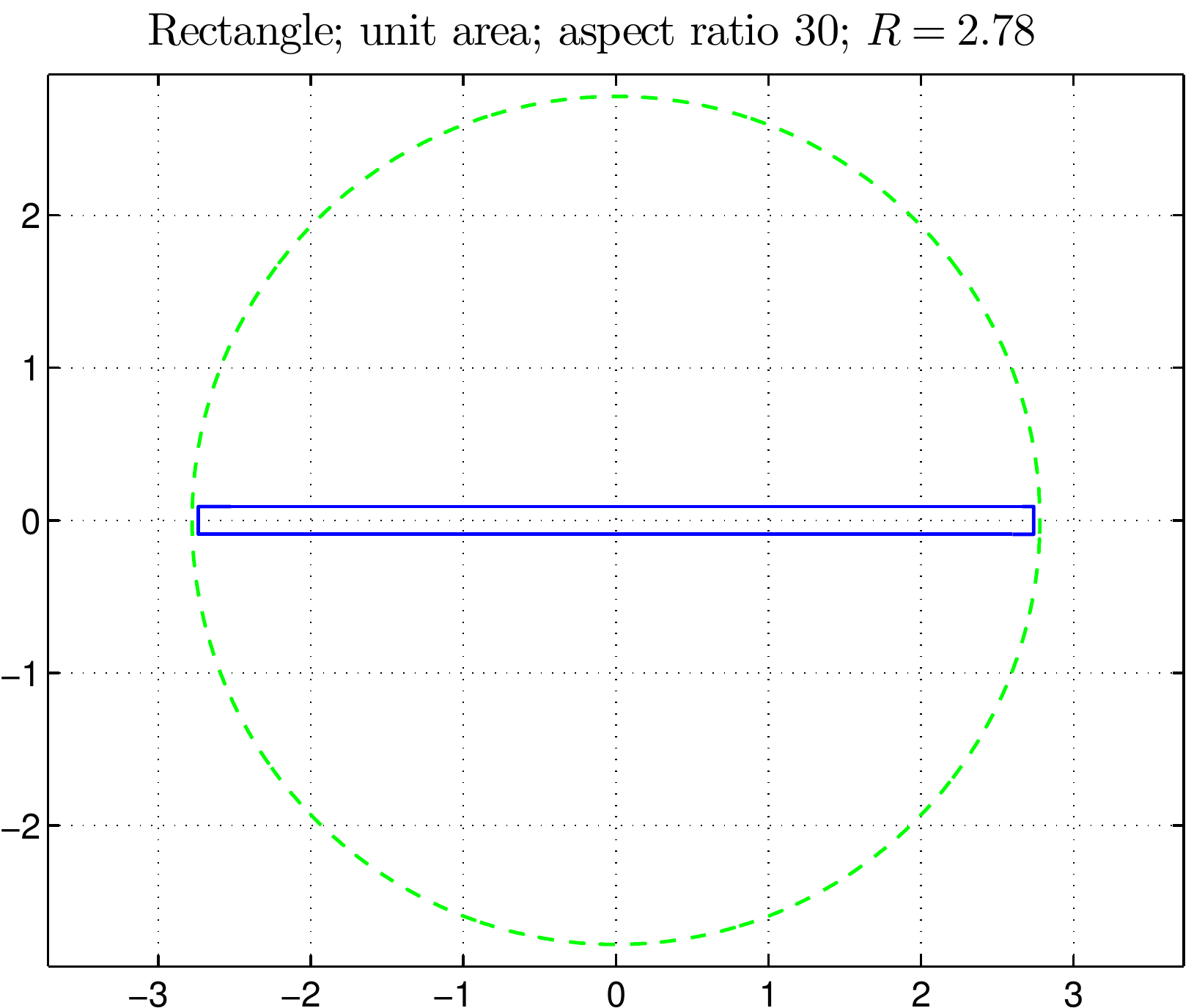}\hskip 7mm
      \includegraphics[width=6.5cm]{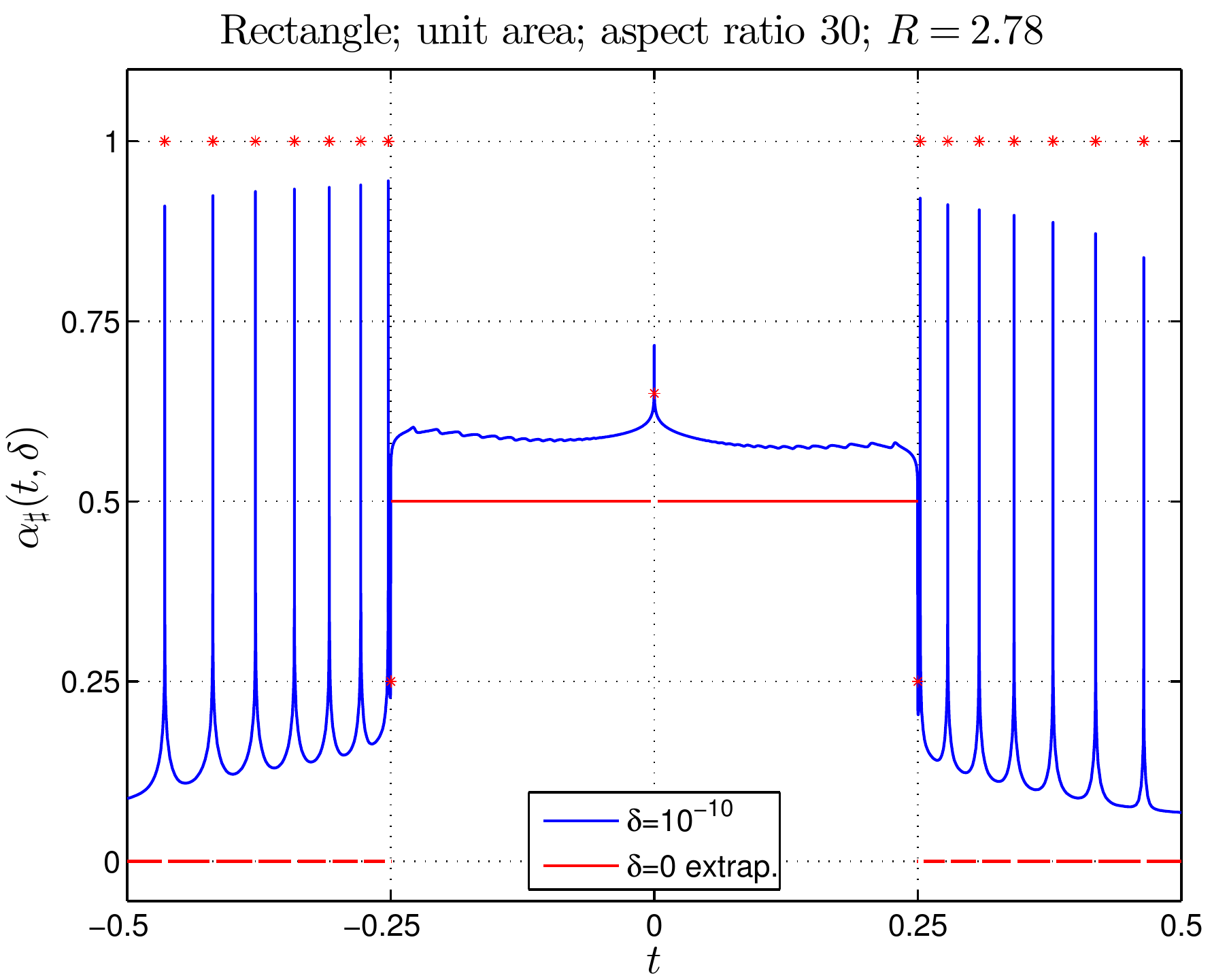}
\end{center}
\caption{Rectangles with various aspect ratios (Left column) and the corresponding spectra (Right column). The second row exhibits a rectangle with the special aspect ratio such that eigenvalues just about to emerge at the two ends of the continuous spectrum interval.}
\label{fig:rectangles}
\end{figure}

The spectrum of the NP operator is computed for rectangles and
superellipses of various aspect ratios $r$. The images of
Fig.\;\ref{fig:rectangles} show $\alpha_\sharp(t)$ for rectangles with
unit area and $r\in\{1; 2.201592; 3; 30\}$. These images illustrate
Corollary~\ref{cor:thinrec}, which says that a rectangle with a
sufficiently high $r$ exhibits eigenvalues outside the continuous
spectrum $t\in[-0.25,0.25]$ and that the number of such eigenvalues
increase with $r$. The ratio $r\approx 2.201592$ is a very special
aspect ratio for which the eigenvalues of the corresponding rectangle
are just about to emerge at $t=\pm 0.25$. It is interesting to observe
that $\alpha_\sharp(t)$ of this rectangle is exactly same as
$\alpha_\sharp(t)$ of the intersecting disks in Fig.\;\ref{fig:ID2} for $t\neq0$.

Superellipses are smooth domains which can be described by the
Cartesian equation
\begin{displaymath}
|x/r|^k+|y|^k=1\,,
\end{displaymath}
where $r$ is the aspect ratio and $k\ge 2$ is a positive parameter.
The higher the parameter $k$ is, the more the superellipse resembles a
rectangle. However, the spectrum of the superellipse always consists
of discrete eigenvalues only. This is so since the corresponding NP
operator is compact. Similarly as with rectangles and ellipses,
superellipses with high $r$ exhibit large eigenvalues.
Table\;\ref{table:analyticRCIP2} shows that eigenvalues of the
superellipse that lie outside the continuous spectrum
$t\in[-0.25,0.25]$ converge to eigenvalues of the rectangle with
the same $r$ as $k\to\infty$. (Eigenvalues of the superellipse
that lie inside the continuous spectrum get increasingly densely
spaced as $k\to\infty$.)

\begin{table}[ht]
\centering 
\begin{tabular}{|c|| l | l |l | l | l |} 
\hline 
$n$&\hskip .8cm$k=10^1$&\hskip .8cm$k=10^4$&\hskip .8cm$k=10^{16}$&\hskip .8cm rectangle\\[0.5ex]
\hline 
1&     0.4641820097578  & 0.4644081276586 & 0.4644081752814 &    0.46440817528139\\[.3ex]
2&  0.4184312731794 &   0.4187549794499 &   0.4187551816213  &   0.41875518162132\\[.3ex]
3&    0.3780806619486  & 0.3783007052834&       0.3783013145614&    0.37830131456136\\ [.3ex] 
4&      0.3413081257441&  0.3413712365784       &0.3413730990324        &0.34137309903240\\[.3ex]
5& 0.3082509222763&0.3082501675778      &  0.3082566649421 &0.30825666494214\\[.3ex]
6&      0.2782621209976 &  0.2783942470929      &  0.278425654617 &0.27842565462101\\[.3ex]
7&  0.2512202243804 &0.2519388130114    &   0.252346607&0.25234907781210\\[.3ex]
8&      0.2267447370526 &  0.2298550809760&      0.247976317&\\
\hline 
\end{tabular}
\caption{Largest eigenvalues $\lambda_n$'s of superellipses $|x/30|^k+|y|^k=1$ and those of rectangle with aspect ratio $30$.}
\label{table:analyticRCIP2}
\end{table}

\subsection{Perturbed ellipse with a corner}\label{subsec:pertell}
\begin{figure}[htbp]
\centering\includegraphics[width=6cm]{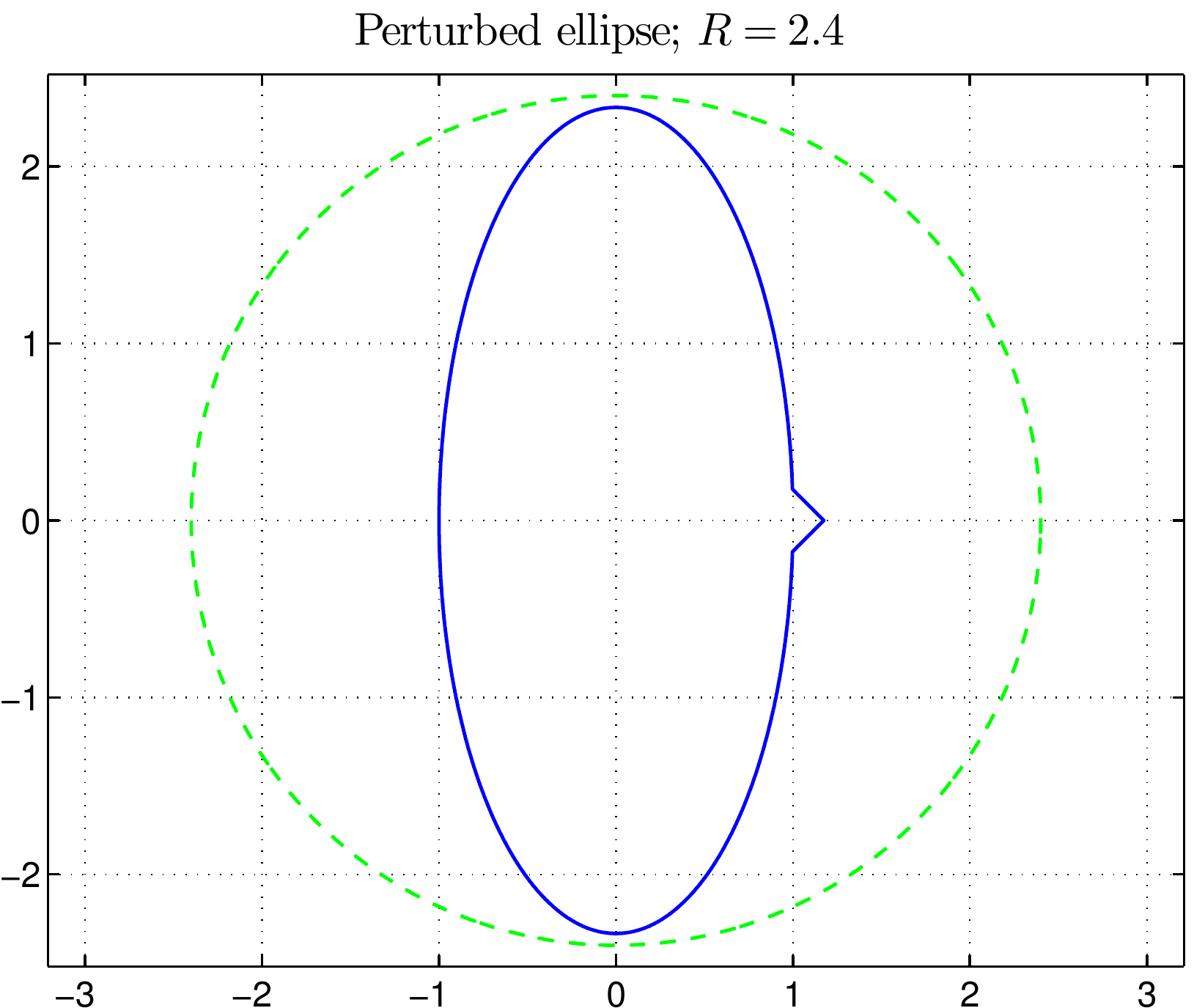}\hskip 7mm
\centering\includegraphics[width=6cm]{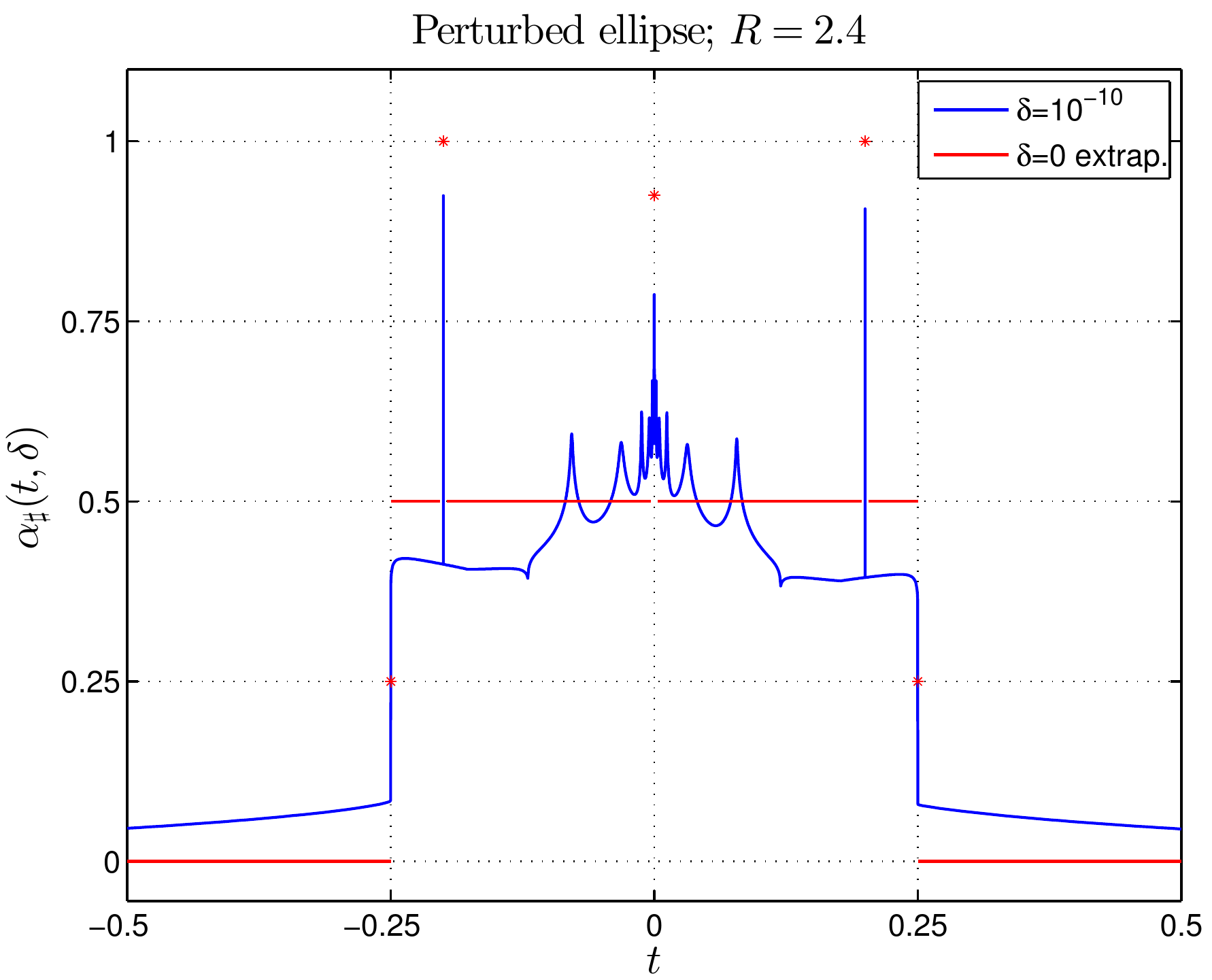}\\[2mm]
\centering\includegraphics[width=6cm]{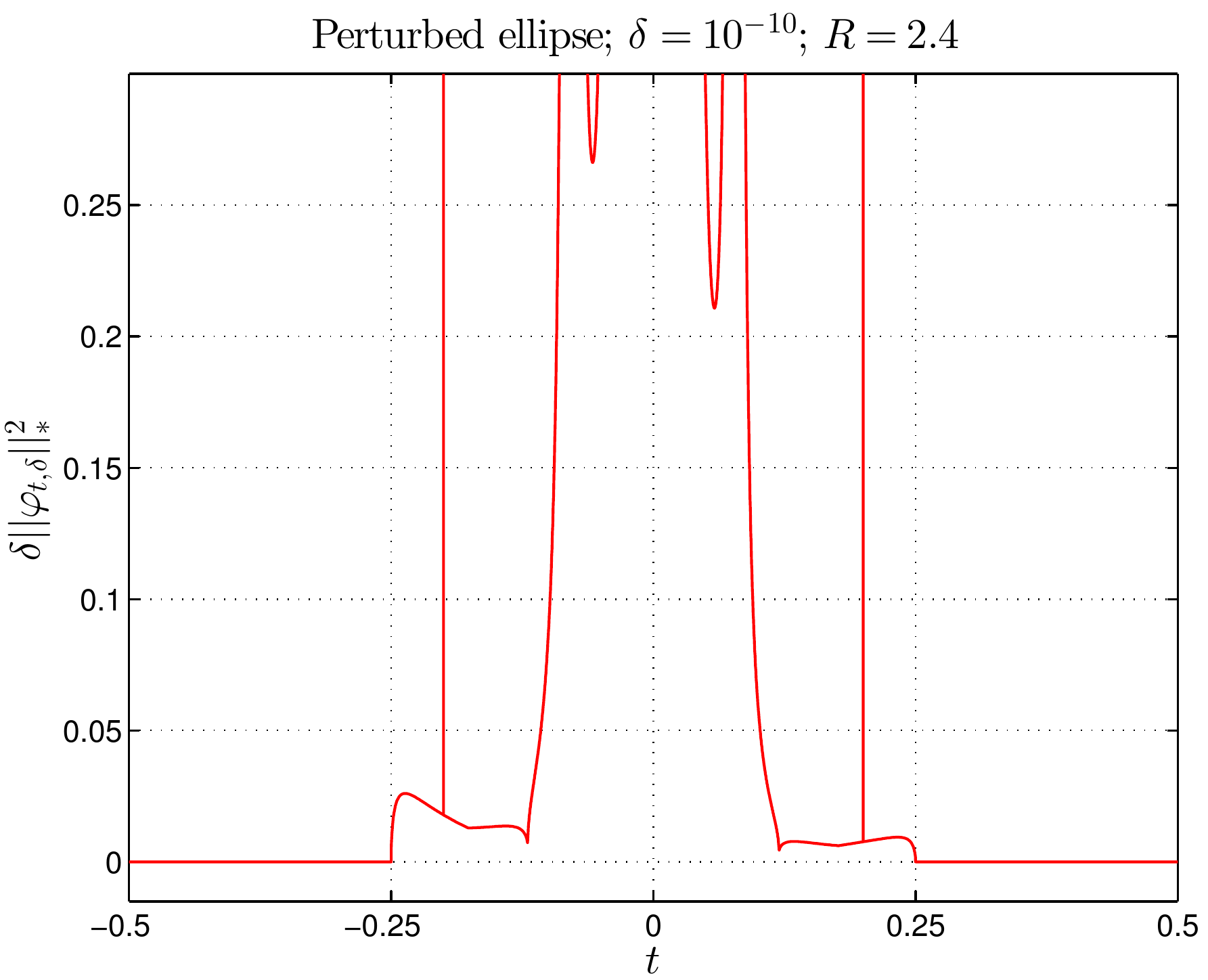}\hskip 7mm
\centering\includegraphics[width=6cm]{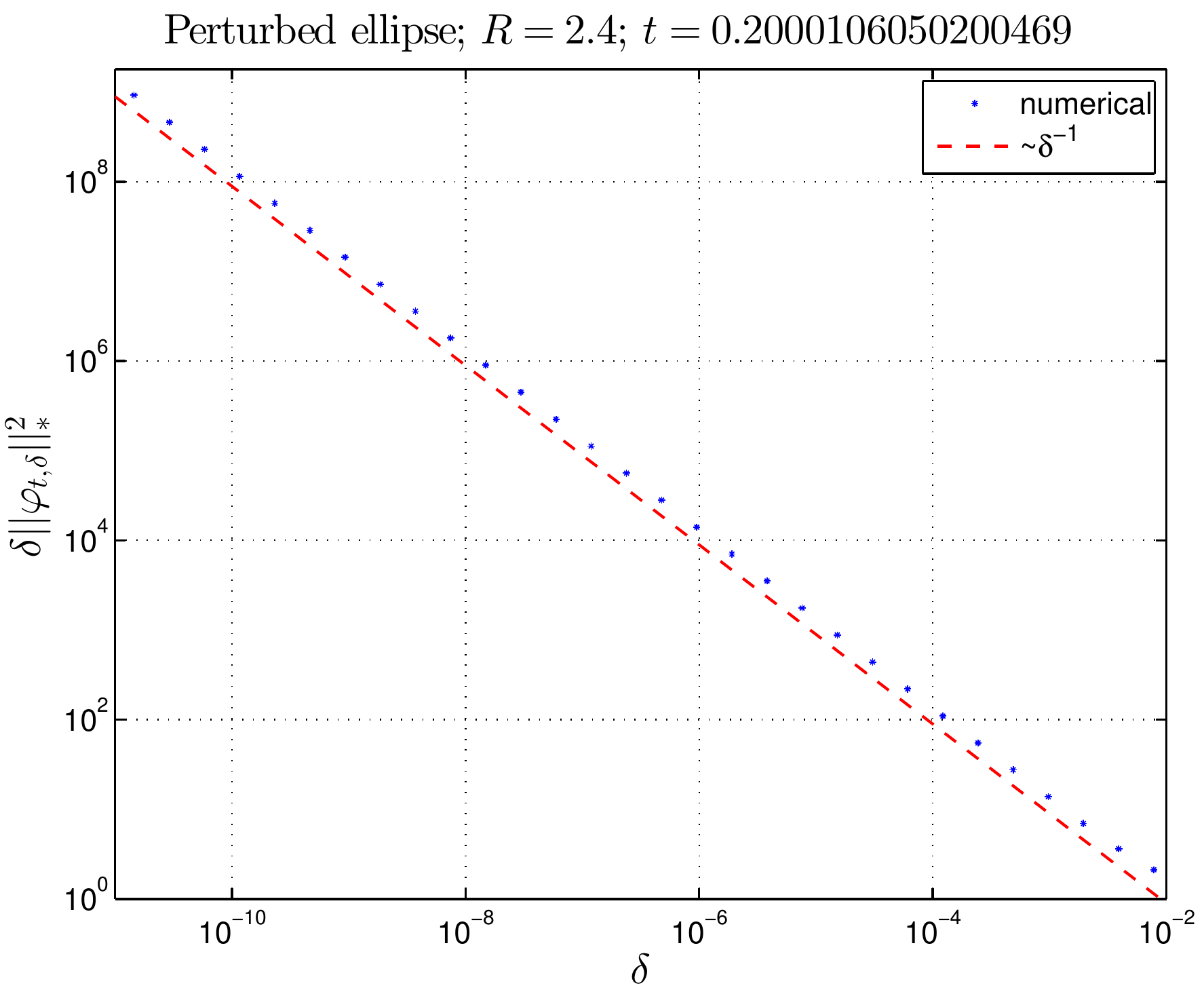}
\caption{Perturbed ellipse with a corner. The aspect ratio before the perturbation is $7/3$. The spectrum contains two singularly continuous spectra.}
\label{fig:PE}
\end{figure}

Fig.\;\ref{fig:PE} shows that even a small Lipschitz perturbation of a smooth domain may induce a big change in the spectrum. The perturbed domain has the interval of continuous spectrum, whose bounds are determined by \eqnref{bounds}. Two singularly continuous spectrum lie inside the continuous spectrum at about $\pm 0.2$, which are the largest eigenvalue of the un-perturbed ellipse.

\section*{Conclusion}

We proposed a method to classify spectra of the NP operator on planar domains with Lipschitz boundaries in terms of resonance rates. The method was implemented  computationally using the RCIP-accelerated Nystr{\"o}m solver on domains such as intersecting disks, a triangle, rectangles, superellipses, and a perturbed ellipse. The results show that the NP operators on all the examples have absolutely continuous spectrum, and some of them have pure point spectrum or singularly continuous spectrum. We also prove rigorously two properties of spectrum suggested by experiments: symmetry of the spectrum and existence of pure point spectrum on rectangles of high aspect ratio.

Several questions are raised by numerical experiments of this paper. On rectangles there is a critical aspect ratio which separates non-existence and existence of eigenvalues, and proving this seems quite interesting. It is also interesting to find the relation between the number of eigenvalues and the aspect ratio. We have shown that the rectangle gets thinner, the spectral bound tends to $1/2$ (Theorem \ref{thm:limbr}). It is interesting and useful to extend this result to general domains. It is desirable to construct in a rigorous manner a domain with corners whose NP operator has singularly continuous spectrum. It is also interesting to show that the NP operator on triangles does not have an eigenvalue.

\section*{Acknowledgement}
We would like to thank Mihai Putinar for sending us the paper \cite{PP-arXiv}, for pointing out existence of Carleman's work \cite{Carleman-book-16}, and for fruitful discussions on the NP operator.


\end{document}